\documentclass[a4paper,Haag duality]{mathscan}
\newtheorem{thm}{Theorem}[section] 
\newtheorem{pro}[thm]{Proposition}  
\newtheorem{cor}[thm]{Corollary}    
\newtheorem{lem}[thm]{Lemma}        
\theoremstyle{definition}           
\newtheorem{rem}[thm]{Remark}       
\newtheorem{defn}[thm]{Definition}  
\newtheorem{exam}[thm]{Example}     

\newcommand{\NI}{\noindent}

\newcommand{\bea}{\begin{eqnarray}}
\newcommand{\eea}{\end{eqnarray}}

\def \b #1 {\bf #1}
\newcommand{\ID}{\mathbb{D}}
\newcommand{\IM}{\mathbb{M}}

\newcommand{\IC}{\mathbb{C}}

\newcommand{\IZ}{\mathbb{Z}}
\newcommand{\cal}{\mathcal}
\newcommand{\cla}{{\cal A}}

\newcommand{\clm}{{\cal M}}
\newcommand{\clz}{{\cal Z}}
\newcommand{\cli}{{\cal I}}
\newcommand{\cls}{{\cal S}}

\newcommand{\clf}{{\cal F}}
\newcommand{\clg}{{\cal G}}
\newcommand{\clh}{{\cal H}}
\newcommand{\clp}{{\cal P}}

\newcommand{\clb}{{\cal B}}

\newcommand{\clj}{{\cal J}}

\newcommand{\clc}{{\cal C}}

\newcommand{\raro}{\rightarrow}

\newcommand{\vsp}{\vskip 1em}

\def \qed {\hfill \vrule height6pt width 6pt depth 0pt}
\newcommand{\be}{\begin{equation}}
\newcommand{\ee}{\end{equation}}
\newcommand{\ben}{\begin{eqnarray*}}
\newcommand{\een}{\end{eqnarray*}}
\begin{document}

\title{ Unital completely positive maps and their operator systems}

\author{ Anilesh Mohari }
\thanks{...}

\address{ The Institute of Mathematical Sciences, CIT Campus, Taramani, Chennai-600113 }

\email{anilesh@imsc.res.in}



\thanks{ I express my gratitude and thanks to Gilles Pisier and \'{E}ric Ricard for pointing out a gap in the proof of Theorem 1.1 given in the first draft with an instructive counter example which helped me to rectify the statement  
to its present form.}  

\begin{abstract}
A vector subspace $\cls$ of $\IM_n(\IC)$ is called unital operator system if $x \in \cls$ if and only if $x^* \in \cls$ and the identity operator $I_n \in \cls$, where $n$ is any fixed 
positive integer. Let $C^*(\cls)$ be the $C^*$ subalgebra of $\IM_n(\IC)$ generated by the operator system $\cls$. We prove that a unital complete order isomorphism $\cli:\cls \raro \cls'$ between two such operator systems $\cls$ and $\cls'$ of $\IM_n(\IC)$ has a unique extension to a $C^*$-isomorphism $\cli:C^*(\cls) \raro C^*(\cls')$ if and only if $\cls$ and $\cls'$ are having equal set of complete ranks. The operator system $\cls = \mbox{span}\{v_iv_j^*:1 \le i,j \le d \}$ is uniquely determined for a unital completely positive map 
$\tau(x)=\sum_{1 \le k \le d} v_kxv_k^*$ of index $d \ge 1$. As an application of our main result, we explore this correspondence and characterize upto cocycle conjugacy all extreme points in the convex set of unital completely positive maps on $\IM_n(\IC)$. Using the main result, we also characterize upto cocycle conjugacy all extreme elements in the convex set 
of normalised trace preserving unital completely positive maps on $\IM_n(\IC)$.       
\end{abstract}

\maketitle 
\section{ Introduction: } 

\vsp 
Let $C(X)$ be the commutative $C^*$-algebra of continuous complex valued functions on a compact Hausdorff space $X$. For two such spaces $X$ and $Y$, we define an endomorphism $\Gamma : C(X) \raro C(Y)$ by
\be 
\Gamma(\psi)(y)= \psi \circ \gamma (y),\;\;y \in Y
\ee
for a continuous map $\gamma:Y \raro X$. In case $\gamma$ is a one to one and onto map, then $\Gamma$ is an automorphism. A well known theorem of M. Stone [Sto2] says that an auto-morphism $\Gamma:C(X) \raro C(Y)$ determines uniquely a continuous one to one and onto map $\gamma:Y \raro X$ such that (1) is valid. A linear subspace $\clf \subseteq C(X)$ is called a {\it function-system} if $\clf$ contains constant functions and $\clf$ is closed under conjugation i.e. $\psi \in \clf$ if and only if $\bar{\psi}(x)=\bar{\psi(x)}$. Another theorem of M. Stone [Sto1] also says that, the closed algebra generated by a function-system $\clf$ is equal to $C(X)$ if and only if $\clf$ separates points of $X$. In particular, this result also proves that an automorphism $\Gamma:C(X) \raro C(Y)$ is determined uniquely by its restriction to $\Gamma:\clf \raro \clg$, where $\clf$ and $\clg$ are two separating-point function-systems for $X$ and $Y$ respectively. 

\vsp 
We fix two function systems $\clf$ and $\clg$ of $C(X)$ and $C(Y)$ respectively. A linear map $\Gamma: \clf \raro \clg$ is called {\it non-negative} if it takes non-negative elements in $\clf$ to non-negative elements in $\clg$. A unital linear map $\Gamma:\clf \raro \clg$ is called an {\it order isomorphism} if it is one to one and onto such that both $\Gamma$ and $\Gamma^{-1}$ are non-negative. If $\Gamma:C(X) \raro C(Y)$ is an automorphism, its restriction on function-systems $\Gamma:\clf \raro \clg$ gives an order isomorphism. Conversely, it is known [LeL] that a unital order isomorphism between two separating function-systems $\clf$ and $\clg$ for $X$ and $Y$ respectively, has an extension to an automorphism between $C(X)$ and $C(Y)$. 
 
\vsp 
In this paper, we address this problem in a more general framework of {\it operator systems}, which we describe in details after briefly developing the mathematical framework of $C^*$-algebras [BRI,Pa].      

\vsp 
A Banach $*$-algebra $\cla$ with norm $||.||$ is called $C^*-$algebra if $||x^*x|| = ||x||^2$. An element $x \in \cla$ is called {\it self-adjoint} if $x^*=x$. An element $x \in \cla$ is called {\it non-negative} if $x=y^*y$ for some $y \in \cla$. We use the notations $\cla_h$ and $\cla_+$ for real vector space of self-adjoint elements in $\cla$ and the convex set of 
non-negative elements of $\cla$ respectively. For two elements $x,y \in \cla_h$, we say 
$x \le y$ if $y-x \ge 0$. We will consider in this paper $C^*$-algebras with units i.e. 
there exists an identity element $I \in \cla$ which satisfies $xI=Ix=x$ for all $x \in \cla$.       

\vsp 
A closed subspace $\clm$ of $\cla$ is called {\it operator space}. A closed vector subspace 
$\cls$ of a unital $C^*$ algebra $\cla$ is called {\it self-adjoint} if $x^* \in \cls$ whenever $x \in \cls$. A self-adjoint subspace $\cls$ of $\cla$ is called {\it operator system} [Ar,Pa] if $I \in \cls$. Let $\cls_+$ be the cone of positive elements in $\cls$ and $C^*(\cls)$ to be the $C^*$-algebra generated by $\cls$. A unital linear map $\cli_0:\cls \raro \cls'$ for two operator systems $\cls \subseteq \cla$ and $\cls' \subseteq \clb$ is called {\it positive} if $\cli_0(\cls_+) \subseteq \cls'_+$. The map $\cli_0: \cls \raro \cls'$ is called {\it completely positive} (CP) if for each $k \ge 1$, $\cli_0^{(k)}=\cli_0 \otimes I_k :\!M_k(\cls) \raro \!M_k(\cls')$, defined by $\cli_0 \otimes I_k [x^i_j]=[ \cli_0(x^i_j) ]$ is positive i.e. $\cli_0^{(k)}(\!M_k(\cls)_+) \subseteq \!M_k(\cls')_+$. Two operator systems $\cls$ and $\cls'$ are called {\it order isomorphic} if there exists $\cli_0:\cls \raro \cls'$, a positive unital one to one and onto linear map such that its inverse is also positive. Two operator systems $\cls$ and $\cls'$ are called {\it complete order isomorphic} if there exists a unital completely positive ( UCP ) $\cli_0:\cls \raro \cls'$ one to one and onto map such that its inverse is also completely positive. We have cited standard text book [Pa] on operator systems for several instants in this text, omitting often details as we have maintained the same  terminology and hopefully also the same notations.  

\vsp 
In a celebrated paper [Ka1] R.V. Kadison proved that an order isomorphism between two $C^*$ algebras is a sum of a {\it morphism} and an {\it anti-morphism} i.e. for two arbitrary $C^*$-algebras $\cla$ and $\clb \subseteq \clb(\clh)$, an order-isomorphism $\cli:\cla \raro \clb$ is a disjoint sum of a morphism and an anti-morphism i.e. there exists a projection $e \in \clb'' \bigcap \clb'$ i.e. center of $\clb$ such that $x \raro \cli(x)e$ is morphism ($*$-linear and multiplicative) and $x \raro \cli(x)(I-e)$ is an anti-morphism ($*$-linear and anti-multiplicative). Thus when $\clb$ is a factor, an order-isomorphism is either an isomorphism or anti-isomorphism. It is a simple observation that anti-morphism part in the decomposition will be absent [ER] if we also demand an order isomorphism to be a CP map. 
For details we refer to Corollary 5.2.3 in [ER]. Thus a complete order isomorphism on $C^*$ algebras is a $C^*$-isomorphism which we call in short isomorphism. In particular, a unital completely positive map ( UCP ) [Fa] $\tau:\!M_n(\IC) \raro \!M_n(\IC)$ is norm preserving i.e. $||\tau(x)||=||x||$ for all $x \in \!M_n(\IC)$ if and only if $\tau(x)=uxu^*$ for some unitary $u \in \!M_n(\IC)$. We also refer to [ChC] for non-unital complete order isomorphism on $C^*$ algebras that need not be $C^*$-isomorphic. Thus throughout this paper, we consider here only unital complete order isomorphism as stated in our definition.   

\vsp 
It is natural to look for a generalization of R.V. Kadison theorem for operator systems in the following sense. Let $\cls \subseteq \cla$ and $\cls' \subseteq \clb$ be two operator systems and $\cli_0: \cls \raro \cls'$ be a unital complete order isomorphism. Is there 
a $C^*$-isomorphism $\cli:C^*(\cls) \raro C^*(\cls')$ extending $\cli_0:\cls \raro \cls'?$
This problem has got its attention ever since William Arveson introduced a Hahn Banach type of extension theorem for completely positive maps on operator systems [Ar1]. A sufficient but not necessary condition is known for quite some time [Ar1]. Isomorphism problem between two operator systems are studied in different contexts with different objectives and motivations [Pa,Pi]. We also wish to draw reader's attention to some interesting related problems reviewed in a recent elegant exposition by F. Douglas on unitary problem [Fa] on similarity. In all these analysis, the isomorphism problem of two specific operator systems is implicitly involved. 

\vsp 
We begin with a negative answer to our problem as follows. Let $S_+$ be the right shift isometry on $\clh_+=l^2(\IZ_+)$ and $\cls$ be the operator system spanned by $\{S_+,S_+^*,I_+ \}$. Let also $\cls'$ be the unitary right shift on $\clh=l^2(\IZ)$ and $\cls'$ be the operator system generated by $\{S,S^*,I \}$. The natural map 
$$\lambda_0 I_+ + \lambda_1 S_+ + \lambda_2 S_+^* \raro \lambda_0 I + \lambda_1 S + \lambda_2 S^*$$ 
is a complete order isomorphism. It is clear that $C^*(\cls)$ and $C^*(\cls')$ are not isomorphic.

\vsp 
Given an operator space $\clm$, a pair $(E,k)$ is called an injective envelope of $\clm$ provided 

\NI (i) $E$ is injective i.e. $E \subseteq \clb(\clh)$ is an injective operator space in 
$\clb(\clh)$ i.e. there exists a completely positive projection $\Phi:\clb(\clh) \raro E$ onto $E$;

\NI (ii) $k:\clm \raro E$ is a complete isometry;

\NI (iii) if $E_1$ is another injective operator space of $\clb(\clh)$ 
with $k(\clm) \subseteq E_1 \subseteq E$, then $E_1=E$. 

\vsp 
We recall M. Hamana basic construction [Ha1,Ha2] ( Theorem 15.4 an 15.6 in [Pa]) which says minimal injective $C^*$ envelope $(E,k)$ for an operator space $\clm$ always exists and unique upto isomorphism. However the algebraic multiplication in $E$ is given by 
$$x \circ y = \Phi(xy)$$

\vsp 
Instead of working with the category of operator spaces, we may as well deal with operator systems and M. Hamana construction ensures existence of a unique upto-isomorphism minimal injective envelop operator system $(k,E)$ of $\cls$ where $k$ is also unital. In the following discussion, we use the notation $C^*_e(\cls)$ for the unique upto isomorphic 
$C^*$ envelope $k:\cls \raro C^*_e(\cls)$ of $\cls$. Thus Theorem 15.6 in [Pa], a unital complete isomorphism $\cli_0:\cls \raro \cls'$ extends to a complete isomorphism $\cli:C^*(\cls) \raro C^*(\cls')$ if $C^*(\cls)$ and $C^*(\cls')$ are minimal injective $C^*$ envelope of $\cls$ and $\cls'$ with respect to their inclusion maps respectively. However, even in the commutative situation $\cls \subseteq C(X)$, $C^*_e(\cls) \neq C^*(\cls)$. For an explicit example, let $\cls$ be the algebra of analytic functions 
on the unit disc $\ID=\{z \in \IC: |z| \le 1 \}$. So $\cls$ is an operator system in $C(\ID)$. The linear map $k: f \raro f_{|} \partial \ID$ is a complete order isomorphism 
of $\cls$ onto $C(\partial \ID)$, where $\partial \ID = \{z \in \IC: |z|=1 \}$. Thus $C^*_e(\cls)= C(\partial \ID) \neq C^*(\cls)$.  

\vsp 
In the following text, we illustrate our problem further in case of operator systems in matrix algebras, We fix an 
orthonormal basis $(e_k)$ for $\IC^3$ and consider the operator system 
$$\cls = \mbox{span}\{ e^1_2(3), e^2_1(3), I_3 \},$$ 
where $e^1_2(3)=|e_1><e_2|$ and $e^2_1(3)=|e_2><e_1|$ are elementary matrices in $\!M_3(\IC)$ 
and $I_3$ is the identity matrix in $\!M_3(\IC)$. So we have $C^*(\cls)= \{ \!M_2(\IC) \oplus z |e_3><e_3| : z \in \IC \}$. We claim that $(i_s,C^*(\cls))$ fails to be 
the minimal injective envelope, where $i_s:\cls \raro C^*(\cls)$ is the inclusion map. For a proof, we consider the operator system 
$$\cls' = \mbox{span} \{ e^1_2(2), e^2_1(2), I_2 \},$$ 
where $e^1_2(2)=|e_1><e_2|,e^2_1(2)=|e_2><e_1|$ are now viewed as elementary matrices 
in $\!M_2(\IC)$. The natural linear map $k:\cls \raro \cls'$ that takes 
$$e^1_2(3),e^2_1(3),I_3 \raro e^1_2(2),e^2_1(2),I_2$$
respectively, is a complete order isomorphism. However $C^*(\cls)$ is not $C^*$ isomorphic 
to $C^*(\cls')=\!M_2(\IC)$. This shows $C^*(\cls')=C^*_e(\cls)$ with $k:\cls \raro \cls'$ 
as its injective map, is as well the $C^*$-envelope of $\cls$. 

\vsp 
Given a unital operator system $\cls \subseteq \IM_n(\IC)$, we will introduce in the next section a family $\{E_k: 1 \le k \le m\}$ of orthogonal projections in the centre of $C^*(\cls)$ with $\sum_{1 \le k \le m} E_k = I_n$ such that $C^*$ envelop of 
$\sum_{j \le k \le m}\cls_k$ is $C^*(\cls_j)$ for $1 \le j \le m$, where $\cls_j=E_j\cls E_j$. The set of integers $(n_1,n_2,..n_m)$, where $(dim(E_k)=n_k,\;1 \le k \le m)$ is called {\it the set of complete ranks } of $\cls$. One of our main results proved in section 2 is stated in the following as theorem.  

\vsp 
\begin{thm}
Let $\cls,\cls'$ be the two unital operator systems of a matrix algebra $\!M_n(\IC)$ and 
$\cli_0:\cls \raro \cls'$ be a unital complete order isomorphism then there exists an isomorphism $\cli:C^*(\cls) \raro C^*(\cls')$ extending $\cli_0:\cls \raro \cls'$ provided $\cls$ and $\cls'$ are having same set of complete ranks.  
\end{thm} 

\vsp 
Theorem 1.1 in particular says that the equality in the set of complete ranks for the operator systems $\cls $ and $\cls'$ is a necessary and sufficient condition for a trace preserving complete order isomorphism $\cli_0:\cls \raro \cls'$ to have a trace preserving $C^*$ isomorphic extension between $C^*(\cls)$ and $C^*(\cls')$. It is not hard to extend Theorem 1.1 to approximately finite (AF) $C^*$-algebra. Theorem 1.1 proved in section 2 is crucial to address our main motivating applications given in section 3 and section 4. Section 3 gives a classification for unital completely positive maps on matrix algebra upto cocycle conjugacy. Section 4 deals with classification problem for unital completely positive trace preserving maps on matrix algebra upto cocycle conjugacy.

\section{Operator systems and their isomorphisms:} 

\vsp 
Given a subset $\cls$ of a unital $C^*$ algebra $\cla$, $C^*(\cls)$ denotes $C^*$-algebra the smallest unital $C^*$-subalgebra 
of $\cla$ containing $\cls$. 

\vsp 
\begin{defn} 
Let $\cls \subseteq \cla$ be an operator system in a unital $C^*$ algebra $\cla$. Given a u.c.p map $\tau:\cls \raro \clb$, we say $\clb$ is a $C^*$-{\it envelope} of $\cls$ and write $\clb=C^*_e(\cls)$ if 

\vsp 
\NI (a) $\tau$ is complete order isomometric;

\vsp 
\NI (b) $\clb=C^*(\tau(\cls))$;

\vsp 
\NI (c) Whenever $\eta:\cls \raro \clc$ is a unital complete isometry with $\clc=C^*(\eta(\cls))$, there exists a $*$-homomorphism $\pi: \clc \raro \clb$ such that $\pi \tau = \eta$ on $\cls$ i.e. 
\end{defn} 

\vsp 
\begin{rem} 
The $*$-homomorphism $\pi:\clc \raro \clb$ is necessarily unital, onto and unique. Since both $\tau$ and $\eta$ are unital, $\pi$ is also unital. Range of $\pi$ is a $C^*$-subalgebra of $\clc$ containing $\tau(\cls)=\pi \eta(\cls)$ and so $\pi(\clc)=\clb$. If $\pi':\clc \raro \clb$ is another $*$-homomorphism satisfying (c) i.e. $\pi' \eta = \pi \eta = \tau$ on $\cls$ then $\pi$ and $\pi'$ are equal on $\eta(\cls)$ and so by multiplicative property of $\pi$ and $\pi'$, they are equal on $C^*(\eta(\cls))$ i.e. on $\clc$.   
\end{rem} 

\vsp 
The $C^*$-envelope is unique upto $*$-isomorphism. 

\vsp 
\begin{pro} 
Let $\cls$ be a unital operator system of a $C^*$-algebra $\cla$ with $C^*$ envelopes 
$\tau:\cls \raro C^*(\tau(\cls))$ and $\eta:\cls \raro C^*(\eta(\cls))$ of $\cls$. Then there exists a $C^*$ isomorphism $\pi:C^*(\tau(\cls)) \raro C^*(\eta(\cls))$ such that $\pi \tau = \eta$ on $\cls$.    
\end{pro} 

\vsp 
\begin{proof} 
There are onto $*$-homomorphisms $\pi:C^*(\tau(\cls)) \raro C^*(\eta(\cls))$ and $\pi':C^*(\eta(s)) \raro C^*(\tau(\cls))$ satisfying $\pi \tau = \eta$ and $\pi' \eta = \tau$ on $\cls$. So $\pi' \pi:C^*(\tau(\cls)) \raro C^*(\tau(\cls))$ is $C^*$ homomorphisms and surjective map with $\pi' \pi \tau = \pi' \eta = \tau$ on $\cls$. Thus by Remark 2.2 on 
the uniqueness of homomorphism, $\pi' \pi $ is the identity map on $C^*(\tau(\cls))$. By the same argument, $\pi \pi'$ is 
the identity map on $C^*(\tau'(\cls))$. Thus $\pi$ is a $*$-isomorphism.     
\end{proof} 

\vsp 
\begin{cor} 
Let $\cls$ be an operator system of a unital $C^*$ algebra $\cla$ and $C^*_e(\cls)$ be its $C^*$-envelope and $\cli :\cls \raro C^*(\cls)$ be the trivial embedding. Then there exists a unique $*$-homomorphism $\pi:C^*(\cls) \raro C^*_e(\cls)$ with range equal to $C^*_e(\cls)$ such that $C^*_e(\cls)=C^*(\cli_{\pi}(\cls))$, where we have used notation $\cli_{\pi}$ for $\pi \cli:\cls \raro C^*_e(\cls)$.    
\end{cor}

\vsp 
Let $J$ be the null space of $\pi$ i.e. $J = \{ x \in C^*(\cls): \pi(x)=0 \}$. So $J$ is a two sided maximal ideal of 
$C^*(\cls)$ and $C^*(\cls) / J $ is isomorphic to $C^*_e(\cls)$. In particular, the restriction of the 
quotient map $q_J:x \raro x + J$ to $\cls$ is a unital order isomorphic map and $\pi q_J = \cli_{\pi}$ on 
$\cls$. 

\vsp 
\begin{defn} 
A two sided ideal $K$ of $\cla=C^*(\cls)$ is called {\it boundary } with respect to $\cls$ if the quotient map $g_K:x \raro x + K$ on $\cla$ is completely isometric once restricted to $\cls$. A maximal boundary ideal of $C^*(\cls)$ is called 
{\it Silov ideal} with respect to $\cls$ in $C^*(\cls)$.  
\end{defn} 

\vsp 
\begin{pro} 
Let $\cls$ be a unital operator system of $\cla=C^*(\cls)$ and $\tau:\cls \raro C^*_e(\cls)$ be a $C^*$-envelope and $\pi:C^*(\cls) \raro C^*_e(\cls)$ be the associated $*$-homomorphism such that $\pi \cli = \tau$ on $\cls$ and $J = \{x \in C^*(\cls):\pi(x)=0 \}$. Then $J$ is a the Silov boundary ideal of $C^*(\cls)$. Conversely, if $J$ is a Silov boundary ideal 
of $C^*(\cls)$ then the quotient map $q_{{J}}: x \raro x+J$ is a C$^*$ envelope i.e. $C^*_e(q_J(\cls))$ is a $C^*$-envelope of $\cls$.  
\end{pro} 

\vsp 
\begin{proof} 
For any boundary ideal $K$ of $C^*(\cls)$, the map $s \raro s+K$ being completely isometric on $\cls$ and so there exists a 
$*$-homomorphism $\pi_K:C^*(\cls / K) \raro C^*_e(\cls)$ with $\pi_K q_k = \tau$. However, $C^*(\cls / K) = C^*(\cls) / K$ and 
$\pi_k(q_K(x)) = \pi(x)$ for all $x \in C^*(\cls)$. So for $x \in K$, we have $\pi(x)= \pi_K(K)=0$. This shows $K \subseteq J$. and thus $J$ is the maximal boundary ideal i.e. Silov boundary. 

\vsp 
Conversely, let $J$ be a maximal boundary ideal of $C^*(\cls)$. So the quotient map $q_J:x \raro x+J$ is completely isometric with a surjective $*$-homomorphism $\pi:C^*(\cls / J) \raro C^*_e(\cls)$ satisfying $\pi q_J=\tau$ on $\cls$. However, 
$C^*(\cls / J) =C^*(\cls) / J$. Let $\sigma:C^*(\cls) \raro C^*_e(\cls)$ be the $*$-homorphism defined by $\sigma(x)=\pi(x+J)$ for $x \in C^*(\cls)$. Then $\sigma = \tau$ on $\cls$ and so $ker(\sigma)$ is a boundary ideal of $C^*(\cls)$ and by maximal property of $J$, we have $ ker(\sigma) \subseteq J $. Thus kernel of $\pi$ is equal to $ker(\sigma) / J$ i.e. trivial. Thus $\pi$ is an isomorphism.  
\end{proof} 

\vsp 
In particular, if $C^*(\cls)$ is simple, then $C^*(\cls)$ is isomorphic to $C^*_e(\cls)$. We have the following simple result.

\begin{cor} 
For two unital operator systems $\cls$ and $\cls'$, if $\cli_0:\cls \raro \cls'$ is a unital order 
isomorphic map and $C^*(\cls)$ and $C^*(\cls')$ are simple then $\cli_0$ has a unique extension to an isomorphic 
map $\hat{\cli}_0:C^*(\cls) \raro C^*(\cls')$.     
\end{cor}

\vsp 
\begin{cor} 
If $\cls \subseteq \IM_n(\IC)$ then there exists a projection $p$ in the centre of $C^*(\cls)$ such that $C^*(\cls_p)$ is a $C^*$-envelope with isometric map $\cli_p:x \raro pxp$ from $\cls$ to $\cls_p$, where $\cls_p=p\cls p$. 
\end{cor}

\vsp 
\begin{proof} 
Let $\tau:\cls \raro C^*_e(\cls)$ be an enveloping algebra of $\cls$ with homomorphism $\pi:C^*(\cls) \raro C^*_e(\cls)$ defined as in Proposition 2.6. So the null space of $\pi$, $J$ being a two sided ideal of $C^*(\cls)$, for a given element $x \in C^*(\cls)$, we may write $x=y+z$ with unique $z \in J$ satisfying $||x-z|| = \inf \{ ||x - z'||: z' \in J  \}$. 
We identify its range with $C^*_e(\cls)$ and write $C^*(\cls) = C^*_e(\cls) \oplus J$ and the map $x \raro y$ is a $*$-homorphism with range isomorphic to $C^*_e(\cls)$. So there exists a unique element $p \in C^*_e(\cls)$ so that $I=p \oplus I-p$ with $I-p \in J$. By uniqueness of the decomposition, $p^*=p$. Using two sided ideal property of $J$ in $C^*(\cls)$ in the identity $p=p^2 \oplus p(I-p)$, $p(I-p)=0$. So $p$ is a projection in $C^*(\cls)$ for which $\pi(p)$ is the identity element in $C^*_e(\cls)$. So the homomorphism $\pi:C^*(\cls) \raro C^*_e(\cls)$ is given by $x \raro pxp$ for some projection $p$ in the centre of $C^*(\cls)$, where we identified $C^*_e(\cls)$ as subalgebra of $C^*(\cls)$ and thus $\tau(x) = pxp$ for all $x \in \cls$. 

\vsp 
Alternatively, since $C^*(\cls)$ is isomorphic to a direct sum of matrix sub-algebras of $\IM_n(\IC)$ and $J$ is a two sided ideal of $C^*(\cls)$, there exists a unique minimal projection $p$ in $C^*(\cls) \bigcap C^*(\cls)'$ for which $\pi(p)$ is the identity in $C^*_e(\cls)$ by semi-simple property of $C^*(\cls)$. Thus $C^*_e(\cls)$ is isomorphic to $C^*(\cls_p)$, where $\cls_p=p\cls p$ and $\tau:\cls \raro \cls_p$ is a complete order isomorphism given by the map that takes $x$ to $pxp$ for 
$x \in \cls$.    
\end{proof} 

\vsp 
For a given unital operator system $\cls$, let $\cls_{\clz}$ be the unique minimal operator system that contains $\cls$ and the centre of $C^*(\cls)$.  

\vsp 
\begin{pro} 
Let $\cls$ and $\cls'$ be two unital operator systems of $\IM_n(\IC)$ and $\cli_0:\cls \raro \cls'$ be a unital complete isometric map. Then $\cli_0$ extends to a complete order isomorphism $\hat{\cli}_0:C^*(\cls) \raro C^*(\cls')$ if and only if $\cli_0$ extends to a complete order isomorphism between $\cls_{\clz}$ and $\cls'_{\clz}$.  
\end{pro}

\vsp 
\begin{proof} We only need to prove `if' part of the statement since the converse statement is obvious. We assume without loss of generality that $\cls =\cls_{\clz}$ and $\cls' = \cls'_{\clz}$. By Corollary 2.7, $C^*_e(\cls)$ and $C^*_e(\cls')$ are isomorphic. 

\vsp 
Let $\pi:C^*(\cls) \raro C^*_e(\cls)$ and $\pi':C^*(\cls') \raro C^*_e(\cls')$ be their associated surjective homomorphisms satisfying $\pi=\tau$ on $\cls$ and $\pi'=\tau'$ on $\cls'$ with $C^*(\cls)=C_e^*(\cls) \oplus \clj$ and $C^*(\cls') = C^*_e(\cls') \oplus J'$ respectively. So we have minimal projections $p \in \cls$ and $p' \in \cls'$ in the centre of $C^*(\cls)$ and $C^*(\cls')$ such that $\pi(p)$ and $\pi'(p')$ are the identity elements 
in $C^*_e(\cls)$ and $C^*_e(\cls')$ with complete order isomorphisms $\tau:\cls \raro \cls_p$ and $\tau':\cls' \raro \cls'_{p'}$ defined by $\tau(x)=pxp,\;\forall x \in \cls$ and $\tau'(y)=p'yp'\;\forall y \in \cls'$ respectively. We identify $C^*_e(\cls)$ and $C^*_e(\cls')$ with sub-algebras $C^*(p \cls p)$ and $C^*(p' \cls' p')$ respectively.

\vsp 
We consider the complete order isomorphism $\cli_{0,p}$ from $\cls$ into $\cls'$ defined by 
$$\cli_{0,p}(x)=\cli_0(pxp),\;x \in \cls$$
So there exists a $*$-homomorphism $\pi_{0,p}:C^*(\cli_{0,p}(\cls)) \raro C^*_e(\cls)$ such 
that $\pi_{0,p} \cli_{0,p}=\tau$ on $\cls$. 
However the map $\pi_{0,p} \cli_{0,p}:\cls_p \subseteq \cls \raro C^*_e(\cls)$ has a unique extension 
to an $*$-automorphism and thus the range of the map $\cli_{0,p}$ is equal to $p' \cls' p'$ for some projection 
$p' \in \cls'$. So there exists a complete order isomorphism $\hat{\cli}_{p,0}:\cls \raro \cls'$ satisfying
\be 
\cli_0(pxp) = p'\hat{\cli}_{0,p}(x)p'
\ee 
for all $x \in \cls$. In particular, $\cli_0(p)=p'$ and $\hat{\cli}_{0,p}(p) \ge p'$. Since $\cli_0(p)\cli_0(p)=\cli_0(p)$, by Kadison-Schwarz inequality, we also have $\cli_0(pxp)=\cli_0(p)\cli_0(xp)=\cli_0(p)\cli_0(x)\cli_0(p)$ for all $x \in \cls$. Since the map $x \raro pxp$ is one to one from $\cls$ to $\cls_p$, we conclude that $\hat{\cli}_0=\cli_0$ on $\cls$. 

\vsp 
The reduced map $(\cli_0)^{p^{\perp}}_{(p')^{\perp}}: p^{\perp}xp^{\perp} \raro \cli_0(p^{\perp}xp^{\perp}) = (p')^{\perp}\cli_0(x)(p')^{\perp})$ for all $x \in \cls$ is a unital complete order-isomorphic map from $\cls_{p^{\perp}}=p^{\perp}\cls p^{\perp}$ onto $\cls'_{(p')^{\perp}} = (p')^{\perp}\cls'(p')^{\perp}$. 

\vsp 
Now we repeat the argument used in the first half to show $C^*_e(\cls_{p^{\perp}})$ and 
$C^*_e(\cls'_{\cli_0(p^{\perp})})$ are isomorphic $C^*$-sub-algebras of $\IM_{n-n_1}(\IC)$. 
Since both $\cls$ and $\cls'$ are operator systems of the same matrix algebra of fixed dimension, 
this process will terminate in finite steps with total orthogonal projections $(p_{\alpha})$ in the centre of $C^*(\cls)$ 
and a $C^*$-isomorphism between $C^*(\cls)$ and $C^*(\cls')$ that takes $p_{\alpha}C^*(\cls)p_{\alpha}$ to $p'_{\alpha}C^*(\cls')p'_{\alpha}$ where $(p'_{\alpha})$ is possibly another total orthogonal projections in the centre of $C^*(\cls')$. 
\end{proof}

\vsp 
So the last proposition raises an interesting non trivial question: when can we expect a complete order isomorphism $\cli_0:\cls \raro \cls'$ between two unital operator systems of $\IM_n(\IC)$ to admit a complete order isometric extension from $\cls_{\clz}$ onto $\cls'_{\clz}$? We begin with an instructive counter example. We consider the operator systems 
$$\cls = \mbox{span} \{ e^1_2(4),e^2_1(4), I_4 \} \;\;\mbox{and} \;\; \cls' =\mbox{span} 
\{e^1_2(4)+e^3_4(4), e^2_1(4)+e^4_3(4), I_4 \}$$ 
in $\!M_4(\IC)$. The linear map that takes 
$$e^1_2(4),e^2_1(4), I_4 \raro  e^1_2(4)+e^3_4(4), e^2_1(4)+e^4_3(4), I_4 $$ 
respectively, is a complete order isomorphism but 
$$C^*(\cls) = \IM_2(\IC) \oplus \IC I_2\;\;\mbox{and}\;\;C^*(\cls')=\IM_2(\IC) \otimes I_2$$ are not isomorphic though their injective envelops are isomorphic to $\IM_2(\IC)$ with linear maps that takes
$$\{e^1_2(4),e^2_1(4),I_4 \} \raro \{e^1_2(2),e^2_1(2),I_2 \}$$
and
$$\{e^1_2(4)+e^2_1(4),e^2_1(4)+e^1_4(4),I_4 \} \raro \{e^1_2(2),e^2_1(2),I_2\}$$
repectively. Though the order isomorphism between $\cls$ and $\cls'$ preserves the 
normalised trace of $\IM_4(\IC)$, centres of $C^*(\cls)$ and $C^*(\cls')$ 
are not isomorphic. In other words, the order isomorphism fails to get an extension between 
$\cls_{\clz}$ and $\cls'_{\clz}$, where we used $\cls_{\clz}$ be the minimal operator system 
that contains $\cls$ and the centre of $C^*(\cls)$.    

\vsp 
For a given unital operator system $\cls$, there exists a minimal projection $p$ in the centre of $C^*(\cls)$ so that the map $x \raro pxp$ is a complete order isomorphism between $\cls$ and $\cls_p$ and $C^*(\cls_p)$ is its $C^*$-envelop for the map $x \raro pxp$ of $\cls$. However, such a choice for a projection $p$ in the centre of $C^*(\cls)$ is unique.
We say dimension of $p$ as rank of $\cls$. We find a family of orthogonal projections $\{p_m:m \ge 1\}$ in the centre of $C^*(\cls_{p^{\perp}_{m-1}})$ such that the map $x \raro p_mxp_m$ is a complete order isomorphism from $\cls_{P^{\perp}_{m-1}} \raro p_m\cls_{P^{\perp}_{m-1}}p_m$, where $P_m=\sum_{1 \le k \le m}p_k$ for $m \ge 1$ with $p_0=0$. In particular, $\sum_m p_m =I_d$, where the sum is over a finite index set. The multi-index 
$(n_m)=(dim(p_m))$ is called the set of {\it complete ranks} of $\cls$.     

\vsp 
\begin{pro} 
Let $\cls$ and $\cls'$ be two unital operator systems of $\IM_n(\IC)$ and $\cli_0:\cls \raro \cls'$ be a unital complete order isomorphism. Then there exists a complete order isomorphic extension of $\cli_0:\cls \raro \cls'$ to $\hat{\cli}_0:\cls_{\clz} \raro \cls_{\clz'}$ if and only if $\cls$ and $\cls'$ are having same set of complete ranks. 
\end{pro}

\begin{proof}
Following the proof for Proposition 2.9, we find projections $p \in C^*(\cls)$ and $p' \in C^*(\cls')$ such that the map $\hat{\cli}_0: p \cls p + \IC (I-p) \raro p' \cls' p' + \IC (I-p')$ defined by 
\be 
\hat{\cli}_0(pxp+cp^{\perp})=p'\cli_0(x)p'+c(p')^{\perp}
\ee
for all $x \in \cls$ and $c \in \IC$. Since the maps $x \raro pxp,\;x \in \cls$ and $y \raro pyp, \; y \in \cls'$ are complete order-isomorphisms from $\cls$ onto $p \cls p$ and $\cls'$ onto $p ' \cls' p'$ respectively, (3) sets a well defined complete order-isomorphism as $p \neq I_n$ if and only if $p' \neq I_n$. So there exists a unital completely positive map $\hat{\cli}_0:\cls_p \raro \IM_n(\IC)$ extending $\hat{\cli}_0:p \cls p +\IC(I-p) \raro p' \cls' p' + \IC (I-p')$ by Arveson's extension theorem [Pa3], where $\cls_p=\{ pxp+(1-p)y(I-p):x,y \in \cls \}$. Since $\hat{\cli}_0(p)=p'$ and $p'$ is a projection in $C^*(\cls')$, by Kadison-Schwarz inequality [Ka2], we have 
$\hat{\cli}_0(xp)=\hat{\cli}_0(x)p'$ for all $x \in \cls_p$ and so 
$$\hat{\cli}_0(pxp+(1-p)y(I-p))$$
$$=\hat{\cli}_0(pxp)+\hat{\cli}_0((I-p)y(I-p))$$
$$=p'\hat{\cli}_0(x)p'+ (1-p')\hat{\cli}_0(y)(I-p')$$
for all $x,y \in \cls$. By taking $y=0$, we get $p'\hat{\cli}_0(x)p'=p'\cli_0(x)p'$ for all $x \in \cls$ and thus $\hat{\cli}_0(x)=p'\cli_0(x)p'+(I-p')\cli_0(x)(I-p')=\cli_0(x)$ 
for all $x \in \cls$. This shows that the map $\hat{\cli}_0:\cls_p \raro \IM_n(\IC)$ is an extension of the map $\cli_0:\cls \raro \cls'$. In particular, it shows that the range of the map $\hat{\cli}_0$ is equal to $\cls'_{p'}$. Now by considering the inverse of the map $\cli_0:\cls \raro \cls'$, we also find a completely positive extension of 
$\cli_0^{-1}:\cls' \raro \cls$ to $\cls'_{p'} \raro \IM_n(\IC)$ so that $\hat{\cli}^{-1}_0 \hat{\cli}_0$ is an extension of the trivial map on $\cls_p$ and thus $\hat{\cli}_0$ is a complete order isomorphism between $\cls_p$ and $\cls'_{p'}$.   

\vsp 
Now we repeat the argument used in Proposition 2.9 with $p_1=p$ and $p_1'=p'$, to find an isomorphism between $C^*(\cls_{p_n})$ and $C^*(\cls'_{p'_n})$ and thus extending the map $\cli_0:\cls \raro \cls'$. This completes the proof.     

\end{proof} 

\vsp 
\begin{rem} 
In (3) we have used $I-p' \neq 0$ if and only if $I-p \neq 0$ for complete order isomorphism property of the map.  
\end{rem}

\vsp 
Now we sum up the main result in the following theorem.

\vsp 
\begin{thm} 
Two unital operator systems $\cls$ and $\cls'$ of equal complete ranks in $\IM_n(\IC)$ are completely order isomorphic if and only if their complete order isomorphism map has isomorphic extension between $C^*(\cls)$ and $C^*(\cls')$. 
\end{thm} 

\vsp 
A subspace $\clm$ of a unital $C^*$-algebra $\cla$ is called operator space. Now we include another application of Theorem 3.4 for operator spaces $\clm$ of $\!M_n(\IC)$ possibly 
without units. For an operator space $\clm$ of a unital $C^*$-algebra, we use notation $C^*(\clm^2)$ for the $C^*$ algebra generated by $\clm^2=\{(x_1x_2...x_n)^*(y_1y_2...y_n): x_i,y_i \in \clm \}$ together with the unit $I$ of $\cla$. We also use the notation $BM(\clm)$ for the bi-module generated by $\clm$ over $C^*(\clm^2)$.   

\vsp 
\begin{thm} 
Let $\clm_1$ and $\clm_2$ be the two operator spaces of $\!M_n(\IC)$ and $\tau_0:\clm_1 \raro \clm_2$ be a bijective completely isometric map and complete ranks of operator systems 
$\clm_1^2$ and $\clm_2^2$ are equal. Then there exists a bijective complete isometric extension $\tau: BM(\clm_1) \raro BM(\clm_2)$ of $\tau_0:\clm_1 \raro \clm_2$ and automorphisms $\alpha_1,\alpha_2:C^*(\clm^2_1) \raro C^*(\clm^2_2)$ satisfying the cocycle relation
$$\tau(abc)=\alpha_1(a)\tau(b)\alpha_2(c)$$
for all $b \in BM(\clm_1)$ and $a,c \in C^*(\clm^2)$
\end{thm}    

\vsp 
\begin{proof} 
We consider the operator systems 
\ben
\cls_k = \left \{ \left( \begin{array}{llll} \lambda &,&\;\; g \;\\ h^* &,&\;\;\mu 
\end{array} \right ),\;\;\lambda, \mu,\in \IC,\; g,h \in \clm_k \right \} \subset \!M_2(\cla_1) 
\een
and the bijective map $\cli_0:\cls_1 \raro \cls_2$, defined by 
$$\cli_0: \left ( \begin{array}{llll} \lambda &,&\;\; g \;\\ h^* &,&\;\;\mu 
\end{array} \right ) \raro \left ( \begin{array}{llll} \lambda &,&\;\; \tau_0(g) \;\\ \tau_0(h)^* &,&\;\;\mu \end{array} \right )$$
Then complete order isomorphism between two unital operator systems 
$\cli_0: \cls_1 \raro \cls_2$ are having equal complete ranks $(2n_1,...,2n_m,..)$ if 
complete ranks of $\clm_1^2$ or $\clm_2^2$ is $(n_1,n_2,..,n_m,..)$. Thus it has a unique extension to a $C^*$ isomorphism $\cli:C^*(\cls_1) \raro C^*(\cls_2)$ by Theorem 2.12. Furthermore, by our construction, we have $\cli(I \otimes e^1_1)=I \otimes e^1_1$ and $\cli(I \otimes e^2_2)=e^2_2$, where $e^1_1$ and $e^2_2$ are canonical projections in $\!M_2(\IC)$ with respect to standard orthonormal basis. So $\cli$ takes corners $C^*(\cls_1)^j_k$ of $C^*(\cls_1)$ to corners $C^*(\cls_2)^j_k$ for each $1 \le j,k \le 2$ respectively. Furthermore, $C^*(\cls_l)^1_1 = C^*(\cls_l)^2_2=C^*(\clm^2)$ and $C^*(\cls_l)^1_2 = (C^*(\cls_l)^2_1)^* = BM(\clm_l)$.

\vsp  
Since $\cli$ takes diagonal block matrices to diagonal block matrices, $\alpha_k:C^*(\clm^2_1) \raro C^*(\clm^2_2),\;k=1,2$ are two automorphisms 
determined uniquely by 
$$\cli(a^1_1 \oplus a^2_2) = \alpha_1(a^1_1) \oplus \alpha_2(a^2_2)$$

\vsp 
Thus the corner map $\tau:BM(\clm_1) \raro BM(\clm_2)$, determined by 
$$\cli: \left ( \begin{array}{llll} \lambda &,&\;\; g \;\\ h^* &,&\;\;\mu 
\end{array} \right ) \raro \left ( \begin{array}{llll} \lambda &,&\;\; \tau(g) \;\\ \tau(h)^* &,&\;\;\mu \end{array} \right ),\;\lambda,\mu \in \IC,\;g,h \in BM(\clm_1)$$
is a bijective isometric extension of the map $\tau_0:\clm_1 \raro \clm_2$. 

\vsp 
Furthermore, that $\tau$ satisfies the following cocycle relation: 
$$\tau(abc)=\alpha_1(a)\tau(b)\alpha_2(c)$$
for $a \in C^*(\clm^2_1),b \in C^*(\clm^2_1),c \in BM(\clm^2_1)$ follows once 
we use multiplicative property of $\cli$ on the element
$$ \left ( \begin{array}{llll} 0 &,&\;\; abc \;\\ 0 &,&\;\; 0 
\end{array} \right ) = \left ( \begin{array}{llll} a &,&\;\; 0 \;\\ 0 &,&\;\; 0 \end{array} \right )
\left ( \begin{array}{llll} 0 &,&\;\; b \;\\ 0 &,&\;\;0 \end{array} \right) 
\left ( \begin{array}{llll} 0 &,&\;\; 0 \;\\ 0 &,&\;\;c \end{array} \right)$$
for $a \in C^*(\clm^2_1), b \in BM(\clm_1), c \in C^*(\clm^2_1).$
\end{proof} 

\vsp 
The unique complete isometric extension $\tau$ of $\tau_0:\clm_1 \raro \clm_2$ is given by $\tau(x)=uxv^*,\; x \in \!M_n(\IC)$ for some unitary matrices $u,v$ in $\!M_n(\IC)$ ( since $\cli:C^*(\cls_1) \raro C^*(\cls_2)$ is a $C^*$-isomorphism on a matrix algebra, it is 
a restriction of an inner automorphism implemented by $u \oplus v$ for some unitary matrices 
$u$ and $v$ on $\IC^n$.)  

\vsp 
We denote by $C^*(\clm)$ the $C^*$-algebra generated by $\clm$ together with $\clm^*=\{x: x^* \in \clm \}$. So $\clm + \clm^*$ is a unital operator system. 

\vsp 
\begin{thm} 
Let $\clm_1$ and $\clm_2$ be two unital operator spaces of $\!M_n(\IC)$ with equal complete ranks of their associated operator systems. If $\cli_0:\clm_1 \raro \clm_2$ is a unital complete isometric map then there exists a unique extension of $\cli_0$ to a $C^*$-isomorphism $\cli:C^*(\clm_1) \raro C^*(\clm_2)$.
\end{thm}

\vsp 
\begin{proof} 
We consider the map $\hat{\cli}_0: \clm_1 + \clm_1^* \raro \clm_2 + \clm_2^*$, defined by 
$\hat{\cli}_0(a + b^*) = \cli_0(a) + \cli_0(b)^*$. By Proposition 2.12 in [Pa], $\hat{\cli}_0$ is a well 
defined positive map, since $\cli_0$ is contractive map. The same argument now also shows 
that the map $\hat{\cli}_0 \otimes I_n$ is a well defined positive map for each $n \ge 0$ since $\cli_0 \otimes I_n$ is contractive. The same argument holds good for the map $\hat{\clj}_0:\clm_2 + \clm^*_2 \raro \clm_1 + \clm^*_1$ defined by 
$$\hat{\clj}_0(a + b^*)=\cli^{-1}_0(a) + \cli^{-1}_0(b)^*$$
for all $a,b \in \clm_2$. Clearly $\hat{\clj}_0$ is the inverse of $\hat{\cli}_0$. 
Thus $\hat{\cli}_0$ extends uniquely to an isomorphism between their $C^*$-algebras by Theorem 2.12.
\end{proof}  

\vsp 
\begin{rem} 
Theorem 2.12 has a ready extension for operator systems $\cls$ and $\cls'$ of $\cla$, where $\cla$ is an AF-$C^*$-algebra.   
\end{rem}

\section{ Extremal elements in $CP_{\sigma}$ and their operator systems:}

\vsp
Let $\clm$ be a von-Neumann algebra acting on a separable Hilbert space $\clh$ over the field of complex numbers. A linear map $\tau: \clm \raro \clm$ is called positive if $\tau(x) \ge 0$ for all $x \ge 0$. Such a map is automatically bounded with norm $||\tau||=||\tau(I)||.$ A map $\tau:\clm \raro \clm$ is called completely positive [St] $(CP)$ if $\tau \otimes I_n: \clm \otimes M_n \raro \clm \otimes M_n$ is positive for each $n \ge 1$ where $\tau \otimes I_n$ is, defined by $(x^i_j) \raro ( \tau(x^i_j) )$ with matrix entries  $(x^i_j)$ are elements in $\clm$. We will use the notation $CP(\clm)$ for the convex set of unital completely positive map. In bounded-weak topology of Arveson [Ar1], $CP(\clm)$ is compact [Pa]. A positive map $\tau:\clm \raro \clm$ is called {\it normal } if $\tau(x_{\alpha}) \uparrow \tau(x)$ whenever $x_{\alpha} \uparrow x$ in weak$^*$-topology of $\clm$. We will use notation $CP_{\sigma}(\clm)$ for the convex set of unital completely positive normal maps on $\clm$. It is a simple observation that $CP_{\sigma}(\clm)$ is a convex face in $CP(\clm)$. An element $\tau \in CP(\clm)$ is called extremal if $\tau=\lambda \tau_1 + (1-\lambda)\tau_0$ for some $\tau_0,\tau_1 \in CP(\clm)$ and $0 < \lambda < 1$ then $\tau_0=\tau_1$. Since $CP_{\sigma}(\clm)$ is a face in $CP(\clm)$, an extremal element in $CP_{\sigma}(\clm)$ is also an extremal element in $CP(\clm)$. 

\vsp 
Two elements $\tau,\tau' \in CP(\clm)$ are said to be cocycle conjugate to each other if 
$\alpha \tau = \tau' \beta$ for some automorphisms $\alpha,\beta$ on $\clm$. It is clear that cocycle conjugacy takes one extremal element to another extremal elements in $CP(\clm)$ and one element of $CP_{\sigma}(\clm)$ to another element of $CP_{\sigma}(\clm)$. In this section we will focus our analysis on $CP_{\sigma}(\clm)$. 
  
\vsp 
A unital normal completely positive map $\tau \in CP_{\sigma}(\clm)$ admits a representation [Ar2,Ch]
\be 
\tau(x)=\sum_k v_k x v_k^*
\ee 
where $\{ v_k:k \ge 1\} $ is a family of contractions on $\clh$ such that $\sum_k v_kv_k^*=1$. However for a given element $\tau$, such a family of operators $\{v^*_k:k \ge 1\}$ is not uniquely determined though the vector space generated by $\{v^*_k:k \ge 1 \}$ over the coefficients in the commutant of $\clm$ is uniquely determined. The dimension of $\clm_{\tau}$, as a vector space with coefficients in $\clm'$, which is independent of the choice of the representation of $\tau$, is called {\it numerical index } or simply the {\it index } of $\tau$. 

\vsp 
In particular, if $\clm$ is the algebra of all bounded operators on a separable Hilbert space $\clh$, then the vector space $\clm_{\tau}$ over $\IC$ generated by $\{ v_k: 1 \le k \le m \}$ is determined uniquely by $\tau$, where could be integer or infinity as well. We can equipe $\clm_{\tau}$ with an inner product $s:\clm_{\tau} \times \clm_{\tau} \raro \IC$ by declairing $s(v^*_i,v_j^*)=\delta^i_j$ provided the family $\{v^*_k:k \ge 1\}$ is linearly independent i.e. for any $(l_k) \in l^2(\IC)$,  we have $l_v=\sum_k l_kv^*_k \in \clm_{\tau}$ and $l_v=0$ if and only if $l=0$. Furthermore, $\{w^*_k, k \ge 1\}$ is another linearly independent family of elements in $\clm_{\tau}$ such that $\tau(x)=\sum_{k \ge 1 }w_kxw_k^*$ if and only if $s(w^*_k,w^*_j)=\delta^k_j$ for all $1 \le j,k$ i.e. there exists an isometric matrix $C=((c^j_k))$ with $C^*C=I$ for $w^*_k=\sum_{j \ge 1}c^k_jv_j^*$. In particular, an element $w^* \in \clm_{\tau}$ if and only if the map 
\be 
x \raro \tau_{\lambda,w}(x)= \tau(x) - \lambda wxw^*
\ee 
is completely positive for some $\lambda > 0$. In such a case, there exists a maximal choice for such a $\lambda > 0$ for which $w \notin \clm_{\tau_{\lambda,w}}$.   

\vsp 
Now we recall criteria given by E. St\o rmer [St\o 1] for an element $\tau$ in $CP_{\sigma}$ to be extremal.

\vsp 
\begin{pro} 
Let $\tau(x)=\sum_{1 \le k \le d} v_ixv_i^*$ be a unital normal $CP$ map on $\clm$. Then 
$\tau$ is extremal in $CP_{\sigma}$ if and only if 
\be 
\sum_{1 \le i,j \le d} v_i\lambda^i_j v^*_j=0 
\ee 
for some $\lambda^i_j \in \clm'$  
if and only if $\lambda^i_j=0$ for all $1 \le i,j \le d$. 
\end{pro} 

\begin{proof}
For a proof, we refer to [St\o 1] and [Ar1]. For a survey [Mo3]. 
\end{proof} 

\vsp 
\begin{pro} Let $\tau$ be an extremal element in $CP_{\sigma}(\clm)$. Then there exists a unique element $\eta \in CP_{\sigma}(\clm)$ so that $\clm_{\eta}=\clm_{\tau}$.
\end{pro} 

\vsp 
\begin{proof} 
Let $\tau(x)=\sum_{ k \ge 1} v_k x v_k^*$ be a representation with element $\{v^*_k: k  \ge 1\}$ in $\clb(\clh)$ linearly independent over $\clm'$ as left multiplication. Let $\eta$ be another element in $CP_{\sigma}$ so that $\clm_{\eta}=\clm_{\tau}$ and we write $\eta(x)= \sum_{1 \le k \le d} l_k x l_k^*$. We choose a matrix $C=(c^i_j)$ with entries in $\clm'$ so that $l^*_k = \sum_j c^k_jv^*_j$ as $\clm_{\eta} = \clm_{\tau}$. Since $\sum v_kv_k^*=1$ and also $\sum_k l_kl_k^*=1$ we get 
$$\sum_{j,j'} v_{j'}( \sum_k \bar{c^k_j} c^{k}_{j'} - \delta^j_{j'} ) v^*_j=0$$ 
Since $\tau$ is an extremal element we get $C^*C = ((\delta^i_j))$ by Proposition 4.1 and thus $\eta=\tau$.   
\end{proof}

\vsp 
Now onwards, we confine our interest to $\clm = \!M_n(\IC)$. 

\vsp 
\begin{thm} 
Let $\tau$ and $\tau'$ be two extremal elements of numerical index $d$ such that operator spaces $\clm_{\tau}$ and $\clm_{\tau'}$ be of equal complete ranks and completely isometric i.e. there exists a completely isometric one to one and onto map $\beta:\clm_{\tau} \raro \clm_{\tau'}$. Then $\tau$ and $\tau'$ are cocycle conjugate i.e. 
$$\tau' \alpha_1(x) = \alpha_2 \tau(x)$$
for all $x \in C^*(\clm^2_{\tau})$ and some automorphisms 
$\alpha_1,\alpha_2:C^*(\clm^2_{\tau}) \raro C^*(\clm^2_{\tau'})$.    
\end{thm} 

\vsp 
\begin{proof} 
By Theorem 2.13, the completely isometric map $\beta:\clm_{\tau} \raro \clm_{\tau'}$ between equal complete ranks operator spaces, has a unique complete isometric extension $\beta:BM(\clm_{\tau}) \raro BM(\clm_{\tau'})$ 
satisfying the cocycle relation 
$$\beta(abc)=\alpha_1(a)\beta(b)\alpha_2(c)$$ 
for all $a,c \in C^*(\clm^2_{\tau})$ and $b \in BM(\clm_{\tau})$, where 
$\alpha_1:C^*(\clm^2_{\tau}) \raro C^*(\clm^2_{\tau'})$ and 
$\alpha_2:C^*(\clm^2_{\tau}) \raro C^*(\clm^2_{\tau'})$ are two automorphisms.   

\vsp 
We will show that $\tau' \alpha_1(x) = \alpha_2 \tau(x)$ for all $x \in C^*(\clm^2_{\tau})$. 
The map $\beta:BM(\clm_{\tau}) \raro BM(\clm_{\tau'})$ be completely isometric, we have 
a complete isometric map $\cli:\clm^2_{\tau} \raro \clm^2_{\tau'}$, defined by 
$$\cli(x^*y)=\beta(x)^*\beta(y)$$
for all $x,y \in \clm_{\tau}$, which has a unique extension to a $*$-isomorphism $\cli:C^*(\clm^2_{\tau}) \raro C^*(\clm^2_{\tau'})$. 

\vsp 
The cocycle property gives the following relations:
$$\beta(b_1c_1)=\beta(b_1)\alpha_2(c_1)$$
So we have 
$$\cli((b_2c_2)^*b_1c_1)$$
$$=\beta(b_2c_2)^*\beta(b_1c_1)$$
$$=\alpha_2(c_2)\beta(b_2)^*\beta(b_1)\alpha_2(c_1)$$
$$=\alpha_2(c_2)^*\cli(b_2^*b_1)\alpha_2(c_1)$$
Since $I \in \cls_{\tau}$, we get
$$\cli(c^*_2c_1)=\alpha_2(c_2)^*\alpha_2(c_1)$$
This shows $\cli=\alpha_2$. Furthermore, we have now
$$\beta(b_1)^* \alpha_1(a) \beta(b_2)$$
$$\beta(b_1)^*\beta(ab_2)$$
$$=\alpha_2(b_1^*a b_2)$$ 
for all $b_1,b_2 \in BM(\clm_{\tau})$ and $a \in C^*(\clm^2_{\tau})$. 

\vsp 
In particular, 
$$\beta(v_i^*) = \beta_g(v'_i)^*= \sum_j g^j_i(v'_j)^*$$
for a matrix $g$ in $\!M_d(\IC)$, we get  
$$\alpha_2(v_iv_j^*)= \beta_g(v'_i)\beta_g(v'_j)^*$$ 
for all $1 \le i,j \le d$. Since $\alpha_2$ is unital, we have $\sum_i \beta_g(v'_i)\beta_g(v'_i)^*=I$ and thus 
$g$ is unitary by Proposition 3.2. 

\vsp 
Now we compute the following identities:
$$\alpha_2(\tau(a)) $$
$$=\sum_k \beta(v_k) \alpha_1(a) \beta(v_k)^*$$
$$= \sum_k \beta_g(v'_k) \alpha_1(a) \beta_g(v'_k)^*$$
$$ = \tau'(\alpha_1(a))$$
for all $a \in C^*(\clm^2_{\tau})$. 
\end{proof} 

\vsp 
Given an element $\tau \in CP(\clm)$ with index $d$, we set operator system 
$$\cls_{\tau}= \{ \sum_{1 \le i,j \le d} \lambda^i_jv_iv_j^*: \lambda^i_j \in \IC \}$$
where $\tau(x)=\sum_{1 \le k \le d}v_kxv_k^*,\;x \in \clm$ is a minimal representation of $\tau$. 
It is clear that $C^*(\cls_{\tau}) \subseteq C^*(\clm^2_{\tau})$. In Theorem 2.13, the complete isometric map $\beta:BM(\clm_{\tau}) \raro BM(\clm_{\tau'})$ canonically gives an automorphism $\alpha_g: C^*(\clm^2_{\tau}) \raro C^*(\clm^2_{\tau'})$ such that $\alpha_g(x^*y)=\beta(x)^*\beta(y)$ for all $x,y \in BM(\clm_{\tau})$. 

\vsp 
By Theorem 2.8, a complete order-isomorphism $\cli_0: \cls_{\tau} \raro \cls_{\tau'}$ between equal complete ranks operator systems, has a unique extension to a $C^*$-isomorphism $\cli: C^*(\cls_{\tau}) \raro C^*(\cls_{\tau'})$. 
Is it necessary that there exists a complete isometric map $\beta:\clm_{\tau} \raro \clm_{\tau'}$ such that 
$$\cli(x^*y)=\beta(x)^*\beta(y)$$
for all $x,y \in \clm_{\tau}$? In other words, does an order-isomorphism betweeon 
two complete ranks of operator systems $\cls_{\tau}$ and $\cls_{\tau'}$, give a cocycle conjugate relation 
between $\tau$ and $\tau'$?  

\vsp 
We aim to classify in the following all extreme points in the convex set of unital completely positive maps on $\!M_2(\IC)$ upto cocycle conjugacy and thus give explicit examples illustrating that a complete order isomorphism between two operator systems 
$\cls_{\tau}$ and $\cls_{\tau'}$ of equal complete ranks, is not a sufficient condition for $\tau$ and $\tau'$ to be cocycle conjugate as shown in following example in the simplest possible case with $d=2$. Thus it demands additional invariance for 
our classification problem.   
 
\vsp 
\begin{exam} 
We will parametrize in the following all extremal elements in $CP(\!M_2(\IC))$ with index $d=2$ upto cocycle conjugacy. For a given element $\tau(x)=v_1xv_1^*+v_2xv_2^*$ with $v_1v_2^*+v_1v_2^*=1$, without loss of generality we may assume that $((tr(v_iv_j^*)))$ is a diagonal matrix i.e. $tr(v_iv^*_j)=\delta^i_j \lambda_i $. Going via a cocycle conjugation, we can assume without loss of generality that $v_1=D_1$ is a diagonal matrix with non-negative entries and $v_2=UD_2$ where $D_1,D_2$ are diagonal matrices with non-negative entries and 
\ben
U= \left ( \begin{array}{llll}&\alpha& ,\;\; \beta \;\; \\  &-\bar{\beta}&,\;\;\bar{\alpha}\;\; \\ 
\end{array} \right )
\een
is an element in $SU(2)$ ( absorbing a phase factor i.e. replacing $v_k$ by $e^{i \theta}v_k$ ), where $|\alpha|^2+|\beta|^2=1$. 

\vsp 
In particular now we get $tr(UD_2D_1)=0$ i.e. 
$$\alpha D_1(1,1)D_2(1,1) + \bar{\alpha} D_1(2,2)D_2(2,2)=0$$ 

\vsp 
\NI Case 1. $\alpha \neq 0$: Taking real and complex part, we get $$D_1(1,1)D_2(1,1)+D_1(2,2)D_2(2,2)=0$$ 
if $Re(\alpha) \neq 0$ and 
$$D_1(1,1)D_2(1,1)-D_1(2,2)D_2(2,2)=0$$ 
if $Im(\alpha) \neq 0$. Unital property also ensures that for $k=1,2$, $$D_k(1,1)^2+D_k(2,2)^2=1$$ 
Now consider the function $f(x)=x(1-x)$ on $[0,1]$ and note that in the later situation $f(D_1(1,1)^2)=f(D_2(1,1)^2)$ and thus $D_1(1,1)=D_2(1,1)$ and thus $D_1=D_2$. By Choi-St\o rmer's criteria Proposition 3.1, $\tau$ is not extremal. Thus we are forced to the situation where all entries are non-negative and 
$$D_1(1,1)D_2(1,1)=D_1(1,1)D_2(2,2)=0$$ 
Thus $U \in SU(2)$ with $Re(\alpha) \neq 0,\;Im(\alpha) = 0$ and either $D_1=|e_1><e_1|,\;D_2=|e_2><e_2|$ or $D_2=|e_2><e_2|,D_1=|e_1><e_1|$. In such a case $v_1=|e_1><e_1|$ and 
$v_2=\beta|e_1><e_2|+\bar{\alpha}|e_2><e_2|$ where $|\alpha|^2+|\beta|^2=1$. But it also shows that $v_2v_1^*=0$. Thus $\tau$ 
is not an extremal element when $\alpha \ne 0$. 

\vsp 
\NI Case 2. $\alpha=0$. Situation is quite simple. It says that $v^*_1=c_1|e_1><e_1|+c_2|e_2><e_2|$ is a pure diagonal and 
$v^*_2 = d_1|e_2><e_1| - d_2 |e_1><e_2|$ is pure off diagonal where $c_1,c_2 \ge 0$ and $d_1,d_2 \ge 0$ and 
where we have used the unitary transformation $e_1 \raro \beta e_1$ and $e_2 \raro e_2$ in order to absorb  
$\beta \in \IC$ with $|\beta|=1$, which does not change the orbit generated by the $CP$ map $\tau$. 

\vsp 
We also compute     
$$v_2v_1^*=UD_2D_1= -c_1d_2 |e_2><e_1| + c_2d_1 |e_1><e_2|$$ 
Thus $\tau$ is an extremal element if these two sets: $v_1v_1^*, v_2v_2^*$ are linearly independent 
and $v_1v_2^*,v_2v_1^*$ are linearly independent. Both give same relation 
$$d_1c_2 \neq d_2c_1$$ as condition for linear independence. 

\vsp 
So orbit space is determined by $\tau(x)=v_1xv_1^*+v_2xv_2^*$ where    
$v_1=c_1|e_1><e_1|+c_2|e_2><e_2|$ 
and 
$v_2=d_1|e_1><e_2| - d_2 |e_2><e_1|$
and parameter space 
satisfies the unital relation   
\be 
c_1^2+d_1^2=1,\;c_2^2 + d_2^2 =1
\ee
with 
$d_1c_2 \neq c_1d_2,\;$
$c_1,c_2,d_1,d_2 \ge 0$

Without loss of generality we assume that 
\be 
c_1^2+c_2^2 \le  1 \le d_1^2+d_2^2
\ee

Further using symmetry without loss of generality we assume that $c_1 < c_2$ and so $d_2 < d_1$. In such a 
case we ensure 
\be 
0 < c_1d_2 -c_2d_1 \neq 0
\ee 

\vsp 
Thus $\tau$ is completely determined by $c_1$ and $c_2$ with range $0 \le c_1 < c_2 \le 1$ where $c_1^2+c_2^2 \le 1$ and each element with distinguished parameters will have non-conjugate orbit since the $((tr(v_iv_j^*))$ is an invariance for cocycle conjugacy but having equal complete rank.    
\end{exam} 

\vsp 
\begin{exam} 
We consider the following elements $v_1,v_2 \in M_4(\IC):$  

\ben
v_1= \left ( \begin{array}{llll} &1&,\; 0,\;\; 0,\;\;0\;\; \\  &0&,\;0,\;\;0,\;\;0\;\;   \\ &0&,\;0,\;{1 \over \sqrt{2}},\;\;0\;\; \\
&0&,\;0,\;\;0,\;\;{1 \over \sqrt{2}}\;\; 
\end{array} \right ),
\een

\ben
v_2= \left (\begin{array}{llll} &0&,\;0,\;\; 0,\;\;0\;\; \\  &0&,\;1,\;\;0,\;\;0\;\;   \\ &0&,\;0,\;{1 \over \sqrt{2}},\;\;0\;\; \\
&0&,\;0,\;\;0,\;\;{i \over \sqrt{2}}\;\; 
\end{array} \right )
\een
So we have $v_1v_1^*=D(1,0,{1 \over 2}, {1 \over 2}),v_1v_2^*=D(0,1,{1 \over 2}, -{i \over 2}),v_2v_1^*=D(0,1,{1 \over 2}, { i \over 2})$ and $v_2v_2^*=D(0,1,{1 \over 2},{1 \over 2 })$. For any given matrix $\lambda=(\lambda^i_j) \in M_2(\IC)$, $\sum_{i,j} \lambda^i_j v_iv_j^* \ge 0$ if and only if 
$\lambda^i_i \ge 0,\; i=1,2,\; \lambda^1_1+\lambda^2_2 +2 Re(\lambda^1_2) \ge 0$ and $\lambda^1_1+\lambda^2_2 + 2Re(i \lambda^1_2) \ge 0$. 
It is simple to show that $\sum \lambda^i_jv_iv_j^* \ge 0$ whenever $\lambda=(\lambda^i_j) \ge 0$. However the above relation says that converse is false. Thus the injective unital map from $M_2(\IC) \raro D_4$, $(4 \times 4)$ diagonal matrices given by $\lambda \raro \sum \lambda^i_jv_iv_j^*$ is not an order isomorphism onto map between two 
operator spaces ( here they are $C^*$-algebras ). Though the map is contractive, inverse map is not so.      

\vsp 
We check that conjugation action $\clj(z_1,z_2,z_3,z_3)=(\bar{z}_1,\bar{z}_2,\bar{z}_3,\bar{z}_4)$ takes $\clj v_1 \clj = v^*_1$ and $\clj v_2 \clj = v_2^*$. The operator spaces generated by the two sets of vectors $\{v_iv_j^*:1 \le i,j \le 2 \}$ and 
$\{v_i^*v_j:1 \le i,j \le 2 \}$ are order isomorphic and they are conjugated by anti-unitary operator $\clj$. 
So $\tau(x)=\sum_{1 \le k \le 2} v_kxv_k^*$ and $\tilde{\tau}(x)=\sum_{1 \le k \le 2}v_k^*xv_k$ are cocycle conjugate by anti-automorphism. 
\end{exam} 

\vsp 
\begin{lem} 
Let $\clj_1,\clj_2$ be two anti-unitary operators on $\clh$. Then the map 
$$\tau \raro \tau^{\clj_1,\clj_1}$$ 
is an affine one to one and onto map on $CP_{\sigma}$ where 
$$\tau^{\clj_1,\clj_2}(x)=\clj_2\tau(\clj_1x\clj_1)\clj_2$$
\end{lem} 

\vsp 
\begin{proof} 
For an anti-unitary operator $\clj$ we have $\clj=\clj^*$ and $\clj^2=I$, where by definition $\clj^*$ is the conjugate linear map, defined by $<\clj^*f,g>= \overline{<f,\clj g>}$, where inner product is taken conjugate linear in the second variable. Thus $$(\clj_1v_k\clj_2)^*=\clj_2v_k^*\clj_1$$ 
So $\tau^{\clj_1,\clj_2}$ is a unital CP map. Other properties are now obvious. 
\end{proof} 

\vsp 
Thus $CP_{\sigma}$ is stable under cocycle conjugation by automorphism as well cocycle conjugation by anti-automorphisms taking extremal elements to extremal elements of $CP_{\sigma}$. 

\vsp 
Let $\tau(x)=\sum_{1 \le k \le d} v_kxv_k^*$ be a minimal representation of an extremal element in $C_{\sigma}(\clm)$. If $d \ge 2$, we can find a matrix $g \in U_d(\IC)$ such that $g (\phi_0(v_kv_j^*)) g^*$ is a diagonal matrix. Thus changing $(v^*_k)$ to $(\beta_g(v_k^*)$, we find a minimal representation $\tau(x)=\sum_{1 \le k \le d}v_kxv_k^*$ with $\phi_0(v_kv_j^*)= \sigma_k \delta^k_j$ for some $0 \le \sigma_k \le 1$ with $\sum_{1 \le k \le d} \sigma_k =1$.     

\vsp 
We begin with the following simple lemma. 

\vsp 
\begin{lem} 
Let $\tau_v(x)=\sum_{1 \le k \le d}v_kxv_k^*$ and $\tau_{v'}(x)=\sum_{1 \le  \le d}(v'_k)x(v'_k)^*$ be the minimal representations of two unital extremal elements in the convex set of unital completely positive maps on $\IM$ in their standard form and $\cls_{\tau}=\cls_{\tau'}$. 
Then there exists a one to one and onto $*$-linear map 
$$\cli_d:(\lambda^k_j) \raro (\lambda')^k_j$$
on $\!M_d(\IC),\;d \ge 1$, determined by the following equality:
\be 
\sum \lambda^k_j v_kv_j^* = \sum (\lambda')^k_jv'_k(v')^*_j
\ee  
Furthermore, $\tau_v$ and $\tau_{v'}$ are cocycle conjugate by automorphisms or anti-automorphisms if and only if the map $\cli_d$ is an order isomorphism on $\!M_d(\IC)$. 
\end{lem}

\begin{proof} 
We assume for the time being index of $\tau$ which is equal to index of $\tau'$ is $d \ge 1$. Since $\!M_d(\IC)$ is a factor, by a theorem of Kadison [Ka1], $\cli$ is either an isomorphism or an anti-isomorphism. So there exists a unitary matrix $g$ on $\IC^d$
satisfying either  
$$\cli(\lambda)=g\lambda g^*,\;\; \lambda \in \!M_d(\IC)$$ 
or 
$$\cli(\lambda)= \clj_0 g  \lambda^* g^*\clj_0, \;\; \lambda \in \!M_d(\IC),$$ 
where $\clj_0$ is the complex conjugation map on $\IC^d$. 
 
\vsp 
In the preceding case, $\cli(\lambda)=g \lambda g^*$, we set the map 
$$\beta_g(v^*_i)f \raro v^*_if$$ 
for all $f \in \clh$ and $1 \le i \le d$ and verify that
the map is inner product preserving since $\beta_g(v_i)\beta_g(v_j^*) = v'_i(v'_j)^*$ for all $1 \le i,j \le d$. Thus the map has a unitary extension $u$ on $\clh$ such that 
$$u\beta_g(v^*_i)=(v'_i)^*$$
This shows $\tau'(x)=\tau(u^*xu)$ for all $x \in \clb(\clh)$.  

\vsp 
In the later case, $\cli(\lambda)=\clj_0 g \lambda^* g^* \clj_0$ is more delicate. We write $\clh = \clh_R \oplus i \clh_R$, where $\clh_R$ is a Hilbert space over the field of real numbers and $\clj_0(f_1+if_2) = f_1-if_2$ is the complex conjugation anti-unitary operator 
on $\clh$. We define the map 
$$v^*_i f \raro \beta_g(v'_i)^*f$$ 
for all $f \in \clh_R$ and $1 \le i \le d$ and verify that the map is anti-inner product preserving i.e. 
$$<\beta_g(v'_i)^* f,\beta_g(v'_j)^* h>$$
$$=<f, v_jv_i^* h>$$
$$=<v_j^*f, v_i^*h>$$
So we may extend the map anti-linearly to an anti-linear operator $\clj:\clh \raro \clh$
such that 
$$\clj v_i^*f = \beta_g(v'_i)^*f$$ 
We verify that $\clj v_i^* \clj_0 = \beta_g(v'_i)^*$ for all $1 \le i \le d $. Thus $\tau$ and $\tau'$ are conjugated by anti-automorphisms. 
\end{proof} 

\vsp 
\begin{thm} 
Let $\tau_v(x)=\sum_i v_ixv_i^*$ and $\tau_{v'}(x)=\sum v_i'x(v_i')^*$ be two extremal elements in $CP_{\sigma}$ with equal numerical index $d \ge 2$ in their standard form and $\cli: \cls_{\tau_{v}} \raro \cls_{\tau_{v'}}$ be a unital complete order isomorphism and the set of complete ranks for $\cls_{\tau_{v}}$ and $\cls_{\tau_{v'}}$ are equal. Then $\tau_v$ and $\tau_{v'}$ are cocycle conjugate either by automorphisms or anti-automorphisms if and only if there exists an order isomorphic map  
$$\cli_d:(\lambda^i_j) \raro (\lambda')^i_j$$
on $\!M_d(\IC)$ determined by the following equality:
\be 
\sum \cli(\lambda^i_j v_iv_j^*) = \sum (\lambda')^i_jv'_i(v')^*_j
\ee  
\end{thm}  

\vsp 
\begin{proof} 
$\cli_d$ being an order isomorphism on $\!M_d(\IC)$, there exists a $g \in U_d(\IC)$ such that 
$$\cli_d(\lambda) = g\lambda g^*\;\;\mbox{or}\; g \lambda^t g^*$$ 
Let $\beta:\clm_{\tau} \raro \clm_{\tau'}$ be the linear map defined by
$$\beta(v^*_i)= \beta_g(v'_i)^*$$   
in the former case, otherwise 
$$\beta(v^*_i)= \beta_{g^t}(v'_i)^*$$
where $\lambda^t,g^t$ are the transposes of $\lambda$ and $g$ with respect 
to the standard orthonormal basis $(e_i)$ of $\IC^d$. Since 
$\cli(x^*y)=\beta(x)^*\beta(y)$ for all $x,y \in \clm_{\tau}$, we verify that 
$\beta$ is a complete isometry as 
$$||[\beta(x^i_j) ]||^2 $$
$$= ||[\beta(x^i_j)]^*[\beta(x^i_j)]|| $$
$$ =||\cli ( [ (x^i_j)]^* [(x^i_j)])|| $$
(by unital complete order-isomorphism) 
$$ =||[(x^i_j)]^*[(x^i_j)]||$$
$$||[(x^i_j)]||^2$$  
Now we complete the proof by Theorem 3.3. 
\end{proof} 

\vsp 
\begin{exam}
(Example 3.4 continued) In this case, $d=2$ and we assume without loss of generality $\tau_v$ and $\tau_{v'}$ are given by the parameters $0 < c_1 < c_2 < 1$ and $0 < c_1' < c_2' < 1$ respectively. So $\tau_v$ and $\tau_{v'}$ are cocycle conjugate to each other if and only if $v'_k = \sum_{1 \le j \le 2}g^k_jv_j$ for some unitary matrix $g=(g^i_j)$. Such an element $g$ exists if and only if $c_1=c'_1$ and $c_2=c'_2$ with $g=I_2$. 
\end{exam}

\vsp 
\begin{exam} 
Let $\tau_v$ and $\tau_{v'}$ be two unital extremal elements in $CP_{\sigma}(\IM_3(\IC))$ on matrix algebra $\IM_3(\IC)$ with equal set of complete ranks $(2,1)$ and indices $d$. So $d$ can not be $1$ or $3$ i.e. $d=2$. So we find two dimensional projections $p$ and $p'$ in $\IM_3(\IC)$ that commutes with elements $\cls_{\tau_v}$ and $\cls_{\tau_{v'}}$ respectively. We write $v_iv_j^*= pv_iv_j^*p \oplus \sigma^i_j$, where $\sigma^i_j = (I-p)v_iv_j^*(I-p) \in \IC$. The matrix $(\sigma^i_j)$ is non-negative definite matrix and without loss of generality we assume $\sigma^i_j=\delta^i_j\sigma_i \ge 0$, where $\sum_i\sigma_i=1$. In case $\sigma_i=0$ for some $1 \le i \le 2$ then $\tau_v$ is a trivial extension of an extreme point described in Example 3.9. Now we consider the case for which $0 < \sigma_i < 1$ and $\sigma^i_j=\delta^i_j\sigma_i$ for $1 \le i,j \le 2$. Without loss of generality we assume $p=|e_1\rangle \langle e_1|+|e_2\rangle \langle e_2|$, $v^*_1=\sum_{1 \le i \le 3}c_i|e_i\rangle \langle e_i|$ and $v^*_2=\sum_{1 \le i \le 3}d_i|Ue_i\rangle \langle e_i|$, where $c_i,d_i \ge 0$ with $c_i^2+d_i^2=1$ and $U$ is an unitary matrix. In case $c_i^2=1$ or $0$, $\tau_v$ is a trivial embedding of an extreme point in $CP_{\sigma}(\IM_2(\IC))$. 

\vsp 
We assume without loss of generality, $0 \le c_i,d_i < 1$ for all $1 \le i \le d$. Since $p$ commutes with $\cls_{\tau_v}$, equivalently commutes with $v_2v_1^*$ and $v_1v_2^*$ ( note that $v_1v_1^*$ and $v_2v_2^*$ are already diagonal matrices). 
Since $v_2v_1^*$ and $v_1v_2^*$ are $U$ modulo pre and post multiplied by diagonal invertible matrices with respect to the 
basis $(e_i)$, the commuting property of $p$ with $\cls_{\tau_v}$ ensures that $p$ commutes with $U$. 

\vsp 
So the vector space generated by $\{pv_iv_j^*p:1 \le i,j \le 2 \}$ is atmost 3 dimensional as subspace of $\IM_2(\IC)$ if $\tau_v$ is not a trivial embedding of an extremal element in $CP_{\sigma}(\IM_2(\IC))$. Now we can follow exactly the method used in the computation in Example 3.4 to conclude 
\be
U= \left ( \begin{array}{llll}&0,& \;\;1, \;\;0 \;\; \\  -&1,&\;\;0,\;\;0\;\; \\ &0,&\;\;0,\;\;1\\ 
\end{array} \right )
\ee 
is an element in $SU(3)$. So we have 
\be 
v^*_1=\sum_{1 \le i \le 3} c_i |e_i\rangle \langle e_i|
\ee 
and 
\be 
v^*_2=U\sum_{1 \le i \le 3} d_i |e_i\rangle \langle e_i|
\ee
where $U$ is given in (12). Thus any extremal element $\tau_{v'}$ with complete ranks $(2,1)$ is 
cocycle conjugate to $\tau_v$ for some $0 < c_i < 1$ and $0 < d_i  < 1$ for which 
$c_i^2+d_i^2=1$ for all $1 \le i \le 3$. We compute
$$
v_1v_1^*= \left ( \begin{array}{llll}&c_1^2,& \;\;0, \;\;\;0 \;\; \\  &0,&\;\;c_2^2,\;\;0\;\; \\ &0,&\;\;0,\;\;c^2_3\\ 
\end{array} \right )
$$

$$
v_2v_2^*= \left ( \begin{array}{llll}&d_1^2,& \;\;0, \;\;\;0 \;\; \\  &0,&\;\;d_2^2,\;\;0\;\; \\ &0,&\;\;0,\;\;d^2_3\\ 
\end{array} \right )
$$

and also

$$
v_1v_2^*= \left ( \begin{array}{llll}&0,& \;\;c_1d_2, \;\;0 \;\; \\  -&c_2d_1&\;\;0,\;\;\;\;\;\;0\;\; \\ &0,&\;\;0,\;\;c_3d_3\\ 
\end{array} \right )
$$

$$
v_2v_1^*= \left ( \begin{array}{llll}0, \;\;&-c_2d_1,& \;\;0 \;\; \\  c_1d_2,\;\;&0,&\;\;\;0\;\; \\ 0,\;\;&0,&\;\;c_3d_3\\ 
\end{array} \right )
$$

\vsp 
So irrespective of the values of $c_3$ and $d_3$, $\tau$ is an extremal element in $CP_{\sigma}(\IM_3(\IC))$ if $c_1d_2 \neq c_2d_1$. In case $c_3d_3=0$ i.e. either $(c_3,d_3)=(1,0)$ or $(0,1)$, then linear indepence of the family is ensured if and only if $c_1d_2 \neq c_2d_1$ and $\tau$ is a trivial extension of an extremal element in $CP_{\sigma}(\IM_2(\IC))$. In case $c_3d_3 \neq 0$ and $c_1d_2=c_2d_1 \neq 0$, then $\tau$ is extremal element if and only ${c_3 \over d_3} \neq {c_1 \over d_1}$.  

\vsp 
We are left to deal with extremal element $\tau_v$ with complete ranks $(3)$. In this case $\tau_v$ 
has to be of index $3$ and so the family $\{v_iv_j^*:1 \le i,j \le 3\}$ spans $\IM_3(\IC)$. We set elements $(l_i)$ in $\IM_3(\IC)$ defined by 
\be
l_1= \left ( \begin{array}{llll}&a_1,& \;\;0, \;\;\;0 \;\; \\  &0,&\;\;a_2,\;\;0\;\; \\ &0,&\;\;0,\;\;\;a_3\\ 
\end{array} \right )
\ee

\be
l_2= \left ( \begin{array}{llll}&0,& \;\;b_1, \;\;0 \;\; \\  &b_2,&\;\;0,\;\;\;\;0\;\; \\ &0,&\;\;0,\;\;\;b_3\\ 
\end{array} \right )
\ee

\be
l_3= \left ( \begin{array}{llll}&0,& \;\;0, \;\;\;\;c_1 \;\; \\  &0,&\;\;c_2,\;\;\;0\;\; \\ &c_3,&\;\;0,\;\;\;\;\;0\;\;\\ 
\end{array} \right )
\ee
where $a_k^2+b_k^2+c_k^2=1$ for all $1 \le k \le 2$ and $0 < a_k,b_k,c_k < 1$. 
The vectors $\{l_il_j^*:0 \le i,j \le 2 \}$ are linearly independent if ${a_1 \over b_1} \neq {a_2 \over b_2}$ or spans $3$ dimensional space if ${a_1 \over b_1}={a_2 \over b_2} \neq {a_3 \over b_3}$. Similar algebraic condition ensures that the family 
$\{l_il_j^*:1 \le i,j \le 3 \}$ spans $\IM_3(\IC)$. So for such parameters $\tau_l(x)=\sum_k l_kxl_k^*$ is an extreme point in $CP_{\sigma}(\IM_3(\IC)$. 

\vsp 
In the following text, we will show that any extremal element $\tau_v$ in $CP_{\sigma}(\IM_3(\IC))$ of complete rank $(3)$ is cocycle conjugate to $\tau_{l}$ for some $l=(l_k)$ given as in (15)-(17). The proof is bit broken into several simplified statements.   

\vsp 
\begin{lem} 
Let $\tau_v(x)=\sum_{1 \le k \le 3}v_kxv_k^*$ be an extremal element in 
$CP_{\sigma}(\IM_3(\IC))$. Then there exists an element $g \in U_3(\IC)$ 
such that $\{ \beta_g(v_k)\beta_g(v_k^*):1 \le k \le 3 \}$ is a commutative family 
of matrices.  
\end{lem}

\vsp 
\begin{proof} 
We fix an orthonormal basis $(e_i)$ for $\IC^3$ and set projections 
$p_{ij}=|e_i\rangle \langle e_i| + |e_j\rangle \langle e_j|$ for all $1 \le i < j \le 3$. 
As a first step, we will prove for some $g \in U_3(\IC)$ that the family $\{p_{ij}\beta_g(v_k)\beta_g(v_k)^*p_{ij}:1 \le k \le 3\}$ is linearly dependent 
for all $i \neq j$. Suppose not. Then for each $g \in U_d(\IC)$, there exists $1 \le i(g) < j(g) \le 3$ such that the family 
$\{p(g)\beta_g(v_k)\beta_g(v_k)^*p(g):1 \le k \le 3 \}$ 
is linearly independent, where $p(g)=p_{i(g)j(g)}$. The map $g \raro p(g)$ is continuous since every infinite sequence has a convergent subsequence as the range of the map has only 
finitely many choices. Since $U_3(\IC)$ is a connected set, we conclude $i(g)=i$ and $j(g)=j$ for some $i < j$. Thus for some $1 \le i < j \le 3$, the family $\{p_{ij}\beta_g(v_k)\beta_g(v_k)^*p_{ij}:1 \le k \le 3\}$ is linearly independent for all $g \in U_3(\IC)$.

\vsp 
We simplify notation and use $p$ for $p_{ij}$ in the following text. The vector space $\clm_p$ generated by $\{pv_jv_k^*p:1 \le j,k \le 3 \}$ is $4$-dimensional. The matrices $\{pv_iv_i^*p:1 \le i \le 3 \}$ being linearly independent, we find a pair $(i_0,j_0)$ with $i_0 \neq j_0$ such that $pv_{i_0}v_{j_0}^*p \neq 0$ and together with $\{pv_iv_i^*p:1 \le i \le 3\}$ spans $\clm_p$. If so then $pv_{j_0}v_{i_0}^*p$ together with $\{pv_iv_i^*p: 1 \le i \le 3 \}$ also spans $\clm_p$. So 
$$pv_{i_0}v_{j_0}^*p= \sum_i \lambda_i pv_iv_i^*p + \lambda pv_{j_0}v^*_{i_0}p$$ 
for some $\lambda_i,\lambda \in \IC$. By taking adjoint on both sides, we also get 
$$pv_{j_0}v_{i_0}^*p = \sum_i \bar{\lambda_i} pv_iv_i^*p + \bar{\lambda}pv_{i_0}v^*_{j_0}p$$  
So we have 
$$(1-|\lambda|^2)v_{i_0}v_{j_0}^*=\sum_{1 \le i \le 3}(\lambda_i \bar{\lambda} +\bar{\lambda}_i)pv_iv_i^*p$$
So we have $|\lambda|=1$ and $\lambda_i \bar{\lambda}+\bar{\lambda}_i=0$ for all $1 \le i \le 3$. We may assume without loss
of generality that $\lambda=1$ otherwise we consider $cv^*_{i_0}$ instead of $v^*_{i_0}$ where $c^2=\lambda$ keeping other elements same. Thus $i(pv_{i_0}v_{j_0}^*p-pv_{j_0}v_{i_0}^*p)$ is a self-adjoint element in the linear span of $\{pv_iv_i^*p:1 \le i \le 3 \}$. The self adjoint element $pv_{i_0}v_{j_0}^*p+pv_{j_0}v_{i_0}^*p$ is not an element in the linear span of 
$\{pv_iv_i^*p,\;1 \le i \le 3 \}$ as $pv_{i_0}v_{j_0}^*p$ is linearly independent of them.  This shows that the family $\{pv_iv_j^*p: i,j \in \{i_0,j_0\} \}$ spans $\clm_p$. For 
$k \neq i_0,j_0$ we write  
$$pv_kv_k^*p = \sum_{i,j \in \{i_0,j_0\}} \lambda^i_jpv_iv_j^*p$$
for a symmetric matrix $(\lambda^i_j)$. We use spectral decomposition of $\lambda$ to write
$$pv_kv_k^*p = \sum_{i=i_0,j_0} \lambda_i p\beta_{g_0}(v_i)\beta_{g_0}(v_i^*)p$$ 
for some unitary matrix $g_0 \in \IM_2(\IC)$. This brings a contradiction since the family 
$\{p\beta_{g}(v_i)\beta_{g}(v_i^*)p: 1 \le i \le 3 \}$ is linearly dependent, where $g$ is the trivial extension of $g_0$. 

\vsp 
Thus for some $g \in U_3(\IC)$ and for each $i \le i < j \le 3$, 
$$\{p_{ij}\beta_g(v_k)\beta_g(v_k)^*p_{ij};1 \le k \le 3\}$$ 
is a linearly dependent commuting family, where the commuting property 
follows, once we use linear dependence of these three elements in 
$$p_{ij}\sum_{1 \le k \le 3}\beta_g(v_k)\beta_g(v_k)^*p_{ij}=p_{ij}$$ 
Since this holds for each $p_{ij}$, we conclude that $\{\beta_g(v_k)\beta_g(v_k)^*:1 \le k \le 3 \}$ is a commutative family of diagonal matrices in the basis $(e_i)$. 
\end{proof} 

\vsp 
Thus $\{\beta_g(v_i)\beta_g(v_i)^*:1 \le i \le 3\}$ are diagonal matrices in the basis $(e_i)$ for some $g \in U_3(\IC)$ and so their linear span is the set of all diagonal matrices in the basis $(e_i)$. We set $l_k^*=\beta_g(v_k)^*$ for $1 \le k \le 3$. Without loss of generality we assume that $l_1$ is diagonal with respect to a standard basis $(e_i)$, otherwise we find a unitary matrix such that $l^*_1=U|l^*_1|$ and reset $(l^*_i)$ to $(U^*l^*_i)$ for the general case. So each $l_kl_k^*$ is diagonal in the representation with respect to the basis $(e_i)$ and commutes with the projection $p_{ij}=|e_i\rangle \langle e_i|+|e_j \rangle \langle e_j|$ for all $1 \le i < j \le 3$.

\vsp
\begin{lem} 
Let $p$ be a two dimensional subspace of $\IC^3$ and 
$$\tau_{v,p}(x) = \sum_{1 \le k \le 3}pv_kxv_k^*p$$
be completely positive map on $\IM_3(\IC)$, where $\tau_{v,p}=p$. Then there exists a unique two dimensional subspace $\clm_{\tau,p}$ of $\clm_{\tau}p$ such that the following hold:

\vsp 
\NI (a) For all $x \in \IM_3(\IC)$ 
$$\tau_{v,p}(x) = \sum_{1 \le k \le 3}pw_kxw_k^*p,$$
where $w_k=\beta_h(v_k)$ for $1 \le k \le 3$ for some $h \in U_3(\IC)$ depending on $p$ and 
$$\sum_{1 \le k \le 2}pw_kw_k^*p=\lambda p\;\mbox{and}\;\;pw_3w_3^*p=(1-\lambda)p$$
for some $0 < \lambda \le 1$. 
\vsp 
\NI (b) The family $\{pw_kw_j^*p:1 \le j,k \le 2 \}$ spans $p\IM_3(\IC)p$;

\vsp 
\NI (c) The family $\{w_kw_k^*:1 \le k \le 3\}$ is commutative and $p \in \{w_kw_k^*:1 \le k \le 3 \}''$;
 
\vsp 
\NI (d) $p \in \{w_k^*:1 \le k \le 2 \}'$ and $(I-p)w^*_3p \neq 0$. 

\end{lem}

\vsp 
\begin{proof} 
We consider the CP map $\tau_{w,p}$ on $\IM_3(\IC)$ defined by
$$\tau_{v,p}(x)=\sum_{1 \le k \le 3}pv_kxv_k^*p, \;x \in \IM_3(\IC)$$
Incase $\tau_{v,p}$ is already an extremal point in the convex set of completely 
positive map $\tau$ on $\IM_3(\IC)$ with $\tau(I)=p$, then its index would be $2$ and 
$\clm_{\tau_{v,p}}=\clm_{\tau_v}p$ and the statement (a) is valid with $\lambda=1$. If the map $\tau_{w,p}$ is not so, then we decompose $\tau_{v,p}$ into a convex combination of two elements of index $2$ and $1$ respectively. The element of index $2$ is not a convex combination of two extremal elements of index $1$ by a variation of Proposition 3.1 as $\{pv_kv_j^*p; 1 \le k,j \le 3 \}$ spans $p\IM_3(\IC)p$. 

\vsp 
A proof for (c) follows along the same line of Lemma 3.11 with a restritive choice in the proof for $g \in U_2(\IC)$ extending trivially as element in $U_3(\IC)$ and then we reset 
$w_k$ to $\beta_g(w_k)$ and $w_k = \beta_g \beta_h(v_k)$ for all $1 \le k \le 3$, where 
$h$ has restrictive choice as stated above.  

\vsp 
For (d), we consider the restriction of the map $\tau_{w,p}$ to $\clm_p=p\IM_3(\IC)p \oplus \IC (I-p)$ and so 
$$\tau_{w,p}=\sum_{1 \le k \le k}pw_kxw_k^*p$$ 
$$=\sum_{1 \le k \le 3}pw_kpxpw_k^*p$$
for all $x \in \clm_p$. Thus $w_k^*p = \sum_{1 \le k \le 2} \lambda^k_j pw_k^*p$ for $1 \le k \le 2$ and $(\lambda^k_j:1 \le k,j \le 2 \}$  are elements in $\clm_{p}'$. The matrices on the right hand side commutes with $p$ and so $p$ commutes with $w_k^*p$, i.e. $pw_k^*p=w_k^*p$ for all $1 \le k \le 2$. 

\vsp 
As in the proof of Lemma 3.11, without loss of generality we assume that $w_1^*$ is also an element in the commutative algebra $\{w_kw^*_k: 1 \le k \le 3 \}''$. So we need to show $p$ also commutes with $w_2$. We have
$I-p)w_2^*p=0$. We write $w_2^*$ in the block form 

\be
w_2^*= \left ( \begin{array}{llll}&l_2^*,&\;\;\;c_2 \;\; \\ &c_3,&\;\;\;d_2\;\;\\
\end{array} \right )
\ee
where $l_2^*=pw_2^*p$ and $d_2=(I-p)w_2^*(I-p)$ and $pw_2^*(1-p) = c_2$ and $(1-p)w_2^*p=c_1=0$. 
Since $w_2w_2^*$ is a diagonal matrixs we get $l_2c_2=0$. We write $l_2^*=u|l_2|$, where $u$ is an extension of the partial isometry that takes $|l^*_2|f$ to $l_2^*f$ for all $f \in \IC^3$ i.e. $u|l_2^*|=l_2^*$ and $l^*_2l_2 = ul_2l^*_2u$. 
Thus we have $ul_2l_2^*uc_2=0$ i.e. $l_2^*uc_2=0$. Since the family $\{l_kl_2^*: 1 \le k \le 2 \}$ is linearly independent, we get $auc_2=0$ for some matrix $a$ of rank $2$. Thus $uc_2=0$ for some unitary matrix $u$ in $\IM_2(\IC)$. This shows $c_2=0$.
Thus $p$ commutes with $w_2^*$ and $w^*_2= l_2^* \oplus d_2(I-p)$.   

\vsp 
If $(1-p)w_3^*p=0$ as well, then $p$ would have commuted with $w_3^*$ as well using the same argument by considering $\{w_iw_j^*: i,j \in \{1,3\}$. Then $p$ would have commuted with all $w_kw_j^*:1 \le j,k \le 3$ and so $p$ would have been a scaler multiple of $I$ by the extremal property of $\tau_w$. This brings a contradiction as $0 < p < I$.  

\end{proof}    

\vsp 
\begin{thm} 
Let $\tau$ be an extremal element in $CP_{\sigma}(\IM_3(\IC))$. Then $\tau$ is cocycle conjugate to $\tau_{l}$, where $l=(l_1,l_2,l_3)$ are defined in (15)-(17). 
\end{thm} 

\vsp 
\begin{proof} 
We apply twice Lemma 3.12 with $p=p_{1j}=|e_1\rangle \langle e_1| + |e_i \rangle \langle e_i|$ for $i=2,3$ to conclude $w_1,w_2$ commutes with $p_{12}$ and $w_1,w_3$ commutes with $p_{13}$. Now we follow Example 3.4 to conclude $w_1=l_1,w_2=l_2$ for some $(a_i)$ and $(b_i)$ given in (15) and (16). Along the same argument we also find $w_3=l_3$ for some $(c_i)$ described as in (17). So $\tau$ is cocycle conjugate to $\tau_l$ for some $l=(l_1,l_2,l_3)$ defined in (15)-(17).   
\end{proof} 
\end{exam} 

\vsp 
\begin{rem} 
It is evident that Theorem 3.13 has a ready generalisation for extremal unital completely positive maps on $\IM_n(\IC)$ of index $n$. Thus it establishes a recursive scheme to find all extremal points in $CP_{\sigma}(\IM_n)$ for any $n \ge 1$.  
\end{rem}

\vsp 
\section{Birkhoff's problem for trace preserving unital Markov maps:}

\vsp 
We start with a brief history of this problem. An $n \times n$ matrix $P=(p(i,j):1 \le i,j \le n)$ is called doubly stochastic if 

\NI (a) all entries $p(i,j) \ge 0$; 

\NI (b) all row sums and column sums are equal to $1$ i.e. $\sum_jp(i,j)=1$ for all $i$ and $\sum_i p(i,j)=1$ for all $j$. 

\vsp 
It is obvious that the set $\cls_n$ of all doubly stochastic matrices forms a compact convex subset in $\!R^{n^2}$. A permutation $\pi$ is a one to one map of the indices $\{1,2,..,n\}$ onto themselves. 
The associated permutation matrix $P_{\pi}$ is defined by $P_{\pi}(i,j)=1$ if $j=\pi(i)$ otherwise $0$. Clearly $P_{\pi} \in \cls_n$. It is also simple to check that $P_{\pi}$ is an extremal element in $\cls_n$. Conversely, D. K\"{o}nig [Ko] and G. Birkhoff [Bi1] proved that an extremal point of $\cls_n$ is a permutation matrix $P_{\pi}$ for some permutation $\pi$ on the set $\{1,2,.,n\}$. By Carath\'{e}dory theorem it follows that all doubly stochastic matrices are convex combination of permutation matrices. However, this representation is not unique in general i.e. $\cls_n$ is not a simplex. G. Birkhoff in his book [Bi2] asked for an extension of this problem to infinite many states.

\vsp 
D. G. Kendall [Ke] settled down this conjecture affirmatively as follows. Let the index set be a countable infinite set and $\cls$ be the convex set of doubly stochastic matrices and $\clp_{\pi}$ be the collection of permutation matrices i.e. matrices having exactly one 
unit element in each row and column, all other entries being equal to zero. $\cls$ is viewed as a subset 
of infinite dimensional matrices $\cls_0$ with entries whose rows and columns have uniform bounded $l^1$ norms: 
$$ \mbox{sup}_{i} \sum_{j} |s(i,j)| < \infty,\;\;\mbox{sup}_{j} \sum_{i}|s(i,j)| < \infty $$      
$\cls_0$ is equipped with the coarsest topology such that the linear maps $l^i_j(s)=s(i,j),\;l^i(s)=\sum_j s(i,j)$ and $l_j(s)=\sum_is(i,j)$ are continuous and thus $\cls_0$ becomes a metrizable topological vector space. D. G. Kendall and J. C. Kiefer proved that $\cls_0$ is equal to closer of the convex combination of permutation matrices, where closer is taken in the coarsest topology described above. 

\vsp 
Within the framework of quantum mechanics for irreversible processes, one major problem is to investigate the same problem extending the scope to stochastic or doubly stochastic maps on a non-commutative algebra of observable namely a $C^*$ algebra or a von-Neumann algebra. In the following text we formulate the problem in a general mathematical framework of $C^*$-algebra or von-Neumann algebras $\clm$. As a first step we investigate this problem when $\clm$ is a matrix algebra over the field of complex numbers. We refer readers to Musat and Haagerup [MuH] and also [Oh] for results on extremal points in the compact convex set $$CP_{\sigma,\phi_0}(\clm) =\{ \tau:\clm \raro \clm, \;\mbox{CP map}, \tau(1)=1, \phi_0 \circ \tau=\phi_0 \}$$ 
where $\clm$ is a finite dimensional matrix algebra over the field of complex numbers i.e. $\clm=\!M_n(\IC)$ and $\phi_0$ is the normalized trace on $\clm$. It is known for quite some time that there are extreme points in $CP_{\sigma,\phi_0}(\clm)$ other then automorphisms 
if $n$ is more than equal to three [LS,KuM]. For $n=3,4$, M. Ohno constructed explicit examples of extremal elements in $CP_{\sigma,\phi_0}(\clm)$ those are not extremal in 
the convex set $CP_{\sigma}(\clm)$ of unital completely positive maps on $\clm=\IM_n(\IC)$. 
However, a complete description or characterization of its extreme points remain unkown.  

\vsp
As a first step, we fix $\clm=\IM_n(\IC)$ and aim to classify extremal elements in $CP_{\sigma,\phi_0}(\clm)$ upto cocycle conjugacy. We begin with recalling 
Landau-Streater's criteria for an element $\tau \in CP_{\sigma,\phi_0}(\clm)$ to be an extremal point in $CP_{\sigma,\phi_0}(\clm)$. 

\vsp 
\begin{pro} 
Let $\tau$ be an element in $CP_{\sigma,\phi_0}(\clm)$ with a minimal Stinespring representation $\tau(x)=\sum_{1 \le k \le d} v_kxv_k^*,\;x \in \clm$ and index $d$. 
Then $\tau$ is an extremal element in $CP_{\sigma,\phi_0}$ if and only if
there exists no non-trivial $\lambda=(\lambda^k_j) \in \IM_d(\IC)$ satisfying the relation 
\be 
\sum_{1 \le j,k \le d} \lambda^k_j v_k v^*_j \oplus v^*_jv_k=0
\ee
\end{pro}

\vsp 
\begin{proof} 
For a proof we refer to original work [LS]. For a more general situation, we refer to a recent publication [Mo3].  
\end{proof} 

\vsp 
For an element $\tau(x)=\sum_{1 \le k \le d}v_kxv_k^*,\; x \in \clm$ in $CP_{\sigma,\phi_0}(\clm)$ with index $d$, we set unital operator system $\cls_{\tau,\tilde{\tau}}$ spanned 
by elements $\{v_iv_j^* \oplus v_j^*v_i:1 \le i,j \le d \}$ in $\clm \oplus \clm$. 
It is also obvious that the operator system is independent of representation that 
we have used for $\tau$: For another representation of $\tau$ with $(l_k)$, we have 
$l^*_k = \sum_{1 \le j \le d}\lambda^k_jv^*_j$ and so $l^*_kl_j \oplus l_jl_k^* \in \cls_{\tau,\tilde{\tau}}$ for all $0 \le k,j \le d$. The following proposition gives some more result. 

\vsp 
\begin{pro} 
Let $\tau$ be an extremal element in $CP_{\sigma,\phi_0}(\clm)$. Then there exists an unique element $\eta \in CP_{\sigma,\phi_0}(\clm)$ so that $\clm_{\eta}=\clm_{\tau}$.
\end{pro} 

\vsp 
\begin{proof} 
Let $\tau(x)=\sum_k v_k x v_k^*$ be an minimal representation. Let $\eta$ be another element in $CP_{\phi}$ so that $\clm_{\eta}=\clm_{\tau}$ and we write $\eta(x)= \sum_{1 \le k \le d} l_k x l_k^*$ and 
$\clm_{\eta}=\clm_{\tau}$. We choose $\lambda=(\lambda^i_j)$ so that $l_k = \sum_j \lambda^k_jv_j$ as $\clm_{\eta} 
= \clm_{\tau}$. Since $\sum_k v^*_k v_k= \sum v_kv_k^*=1$ and also $\sum_k l_k^*l_k = \sum_k l_kl_k^*=1$ we get 

$$\sum_{j,j'} ( \sum_k \bar{\lambda^k_j} \lambda^{k}_{j'} - \delta^j_{j'} ) v^*_jv_{j'}=0$$
$$\sum_{j,j'} ( \sum_k \bar{\lambda^k_j} \lambda^{k}_{j'} - \delta^j_{j'} ) v_{j'}v^*_j=0$$ 

Since $\tau$ is an extremal element we get $\lambda \in U_d(\!C)$ by Proposition 4.1 and thus $\eta=\tau$. 
\end{proof}  

\vsp 
\begin{cor} 
Let $\tau$ be an extremal element in $CP_{\sigma,\phi_0}(\clm)$. Then an element $\eta \in CP_{\sigma,\phi_0}(\clm)$ is cocycle conjugate to $\tau$ unitarily or anti-unitarily if and only if $\clm_{\eta}$ is cocycle conjugate to $\clm_{\tau}$  unitarily or anti-unitarily. 
\end{cor} 

\vsp 
\begin{thm} 
Let $\tau(x)=\sum_{1 \le k \le d}v_kxv_k^*$ and $\eta(x)=\sum_{1 \le k \le d} l_kxl_k^*$ 
for $x \in \clm$ be two extremal elements in $CP_{\sigma,\phi_0}(\clm)$ of index $d$. Then $\tau$ and $\eta$ are cocycle conjugate unitarily if and only if unital operator systems 
$\cls_{\tau}$ and $\cls_{\tilde{\tau}}$ are having the set of complete ranks equal to that 
of $\cls_{\eta}$ and $\cls_{\tilde{\eta}}$ respectively and the maps 
$$\cli_g: v_iv_j^* \raro \beta_g(l_i)\beta_g(l_j)^*$$ 
and 
$$\cli_g: v_j^*v_i \raro \beta_g(l_j)^*\beta_g(l_i)$$ 
extends to a complete order-isomorphism from $\cls_{\tau,\tilde{\tau}}$ onto $\cls_{\eta,\tilde{\eta}}$ for some $g \in U_d(\IC)$. Similar statement holds also for anti-unitarily cocycle conjugate extremal elements $\tau$ and $\eta$ in $CP_{\sigma,\phi_0}(\clm)$.   
\end{thm} 

\vsp 
\begin{proof} 
We use Theorem 2.12 with $\cls=\cls_{\tau,\tilde{\tau}}$ and $\cls'=\cls_{\eta,\tilde{\eta}}$ to find a unitary matrix $U=u \oplus v$ in $\clm \oplus \clm$ for $\hat{\cli}_g(X)=UXU^*$ for all $X \in C^*(\cls)$.    
\end{proof} 

\vsp 
In Theorem 4.3, without loss of generality, we may assume that representations of $\tau$ and $\eta$ are in standard form i.e. 
matrices $(\phi_0(v_kv_j^*))$ and $((\phi_0(l_kl_j^*)))$ are diagonal matrices. So $g$ is a unitary matrix satisfing  
$g ((\phi_0(v_kv_j^*))) g^* = (( \phi_0(l_kl_j^*) ))$ i.e. these two diagonal matrices with strict positive values are having 
equal eigen values. So without loss of generality, we may assume that their associated covariance matrices are equal. Furthermore, if the diagonal matrices are having distinct eigen values then $g$ is as well a diagonal matrix modulo a permutation matrix. In such a case, without loss of generality, we may take $g=I_d$ in Theorem 4.3.   

\vsp
\begin{exam} 
Let $\tau$ be an extremal element of $CP_{\sigma,\phi_0}$ in its standard form with index $d$ and $\clm=\IM_3(\IC)$. It is easy to check that $d^2 \le 18$. So $d$ can be at most $4$. Since $\tau \in CP_{\sigma}(\clm)$, we get an extremal decomposition of $\tau$ in $CP_{\sigma}(\clm)$ say $\tau = \sum_k \lambda_k \tau_k$, where each $\tau_k$ is an extremal element in $CP_{\sigma}(\clm)$. 

\vsp 
\NI (a) For $d=2$, $\tau$ is itself an extremal element in $CP_{\sigma}$ as otherwise $\tau$ is a proper convex combination of two automorphisms and thus not extremal in $CP_{\sigma,\phi_0}$. So $\tau$ is cocycle conjugate to $\tau_v(x)=v_1xv_1^*+v_2xv_2^*$ defined in (13)
in Example 3.10 with additional constraint $\tilde{\tau_v}(I)=I$ i.e. $v^*_1v_1+v_2^*v_2 =I$ i.e. $c_1^2+d_2^2=c_2^2+d_1^2=1$. Since $c_i^2 +d_i^2=1$ for all $1 \le i \le 3$, only possible solution is $c_1=c_2=c$ and $d_1=d_2=d$. We assume without loss of generality that 
$0 < c < d < 1$ and compute the following four matrices 
$$
v_1v_1^*= \left ( \begin{array}{llll}&c^2,& \;\;0, \;\;\;0 \;\; \\  &0,&\;\;c^2,\;\;0\;\; \\ &0,&\;\;0,\;\;c^2_3\\ 
\end{array} \right )
$$

$$
v_2v_2^*= \left ( \begin{array}{llll}&d^2,& \;\;0, \;\;\;\;0 \;\; \\  &0,&\;\;d^2,\;\;0\;\; \\ &0,&\;\;0,\;\;d^2_3\\ 
\end{array} \right )
$$

$$
v_1v_2^*= \left ( \begin{array}{llll}&0,& \;\;cd, \;\;0 \;\; \\  -&cd&\;\;0,\;\;\;\;0\;\; \\ &0,&\;\;0,\;\;c_3d_3\\ 
\end{array} \right )
$$

$$
v_2v_1^*= \left ( \begin{array}{llll}&0,& \;\;-cd, \;\;0 \;\; \\  &cd,&\;\;\;\;0,\;\;\;\;0\;\; \\ &0,&\;\;\;\;0,\;\;c_3d_3\\ 
\end{array} \right )
$$
Thus these matrices are linearly independent if $c_3d_3 \neq 0$ and $cd_3 \neq c_3d$. So any extremal element $\tau'$ of index $2$ is cocycle conjugate to $\tau_v$, where $v_1,v_2$ are defined in (13) and (14) with $c_1=c_2=c$ and $d_1=d_2=d$ satisfying $c_3d_3 \neq 0$, $cd_3 \neq c_3d$ and $c_i^2+d_i^2=1$ for all $1 \le i \le 3$.

\vsp 
\NI (b) For $d=3$, if $\tau$ is not an extremal element in $CP_{\sigma}(\clm)$, we write 
$\tau = \lambda \tau_1 + (1-\lambda)\tau_0$ for two elements $\tau_0$ and $\tau_1$ in 
$CP_{\sigma}(\IM_3(\IC))$ and $0 < \lambda < 1$. The vector space $\clm_{\tau}$ is spanned 
by $\clm_{\tau_0}$ and $\clm_{\tau_1}$. Thus indices of $\tau_1$ and $\tau_0$ are less than equal to $2$. So $\tau_1$ and $\tau_0$ are extremal elements in $CP_{\sigma}(\IM_3(\IC))$ of indics $2$, otherwise, say $\tau_0$ is not an extremal element, then $\tau_0$ is at most a convex sum of two automorphisms if not itself is an automorphism and so $\tau_0 \in C_{\sigma,\phi_0}$. Then we have 
$$\phi_0 =\phi_0 \tau$$
$$=\lambda \phi_0 \tau_1 + (1 -\lambda)\phi_0 \tau_0 $$
$$= \lambda \phi_0 \tau_1 + (1-\lambda)\phi_0$$
i.e. $\phi_0 \tau_1 = \phi_0$ as $0 < \lambda < 1$, which contradicts extremal property of 
$\tau$ in $CP_{\sigma,\phi_0}(\IM_3(\IC))$. So we may write
$\tau_0(x) = w_1xw_1 + w_2xw_2^*$ and $\tau_1(x)=w_2xw_2^*+w_3xw_3^*$ 
for some $w_1,w_2 \in \clm_{\tau_0}$ and $w_3,w_4 \in \clm_{\tau_1}$. But $\clm_{\tau_0} \bigcap \clm_{\tau_1}$ is one-dimensional. Thus there exits $l_2=c_1w_1+c_2w_2=c_3w_3+c_4w_4$
for some constants $c_k \in \IC,\;1 \le k \le 4$. We find elements $l_1 \in \clm_{\tau_0}$, $l_3 \in \clm_{\tau_1}$ and constants $a_0,a_1 > 0$ such that 
\be 
\tau_0(x) = l_1xl_1^*+a_0 l_2xl_2^*
\ee 
and 
\be 
\tau_1(x) = a_1 l_2xl_2^* + l_3xl_3^*
\ee
for all $x \in \IM_3(\IC)$. So 
$$\tau(x)=\lambda l_1xl_1^* + (\lambda a_0 + (1-\lambda)a_1)l_2xl_2^* + (1-\lambda)l_3xl_3^*$$
for all $x \in \IM_2(\IC)$. However, both the vector spaces generated by $\{l_kl_k^*:1 \le k \le 2\}$ and $\{l_kl_k^*: 2 \le k \le 3\}$ are two dimensional with their intersection containing atleast $\{I, l_2l_2^*\}$. If $l_2l_2^* = c_2 I$ for some $c_2 > 0$, then $\tau$ is not an extremal element in $CP_{\sigma,\phi_0}$. So we are left to consider the 
situation when $l_2l_2^*$ and $I$ are linearly independent. In such a case, linear spans 
$\{l_kl_k^*:1 \le k \le 2 \}$ and $\{l_kl_k^*: 2 \le k \le 3 \}$ are equal having dimension $2$. Thus $\{l_1l_1^*,l_3l_3^*\}$ are linearly dependent say $l_3l_3^*=cl_1l_1^*$. But $a_1l_1l_1^*-a_0l_3l_3^*=(a_1-a_0)I$ 
and so 
$(a_1-a_0c_1)l_1l_1^*=(a_1-a_0)I.$ 
So $l_1l_1^*=c_1I$ for some $c_1 > 0$ unless $a_0=a_1$. That brings a contradiction to 
extremal property of $\tau$ in $CP_{\sigma,\phi_0}(\IM_3(\IC))$. So $\tau$ is as well 
an extremal element in $CP_{\sigma}(\IM_3(\IC))$.  

\vsp 
Thus in (19) and (20), we may assume $a_0=a_1$. Since $\lambda \tau_1+(1-\lambda)\tau_0=\tau \in CP_{\phi}(\IM_3(\IC))$, we also have
\be 
\lambda l_1^*l_1 + l_2^*l_2 + (1-\lambda)l_3^*l_3=I
\ee
Without loss of generality, we fix a basis $(e_i)$ for $\IC^3$ and assume $l_2^*$ is diagonal
in the basis $(e_i)$ as in Example 3.10. Thus $l_2l_2^*=l^*_2l_2$. We write polar decompostions $l_1^* =U_1|X|$ and $l_3^*=U_3|X|$, where $X$ is the unique positive root of $l_1l^*_1=l_3l_3^*=I-l_2l_2^*=X$. We may assume without loss of generality that $U_1$ and $U_3$ are unitary matrices and rewrite (21) as
\be 
\lambda U_1XU^*_1+(1-\lambda)U_3XU^*_3 = X
\ee 
Note that the map $x \raro \lambda U_1xU^*_1+(1-\lambda)U_3xU_3^*$ is unital and 
tracial state $\phi_0$ preserving. Thus by a standard result in non-commutative ergodic theory [Mo1,Mo2], $X \in \{U_1,U_2 \}'$ and so $l_1l_1^*=l_1^*l_1$ and $l_3l_3^*=l^*_3l_3$. 
However, $\tau$ being an extremal element in $CP_{\sigma,\phi_0}(\IM_3(\IC))$, the family $\{l_kl_k^* \oplus l^*_kl_k: 1 \le k \le 3 \}$ needs to be linearly independent. This contracts equality $l_1l_1=l_3l_3^*$. 

\vsp 
So we conclude that $\tau$ is an extremal element in $CP_{\sigma}(\IM_3(\IC))$. So we have shown any extremal element $\tau$ in $CP_{\sigma,\phi_0}(\IM_3(\IC))$ of index $1 \le 
d \le 3$ is also an extremal element in $CP_{\sigma}(\IM_3(\IC))$. So it shows why H. Ohno
[Oh] had to look for an extremal element $\tau$ in $CP_{\sigma,\phi_0}(\IM_3(\IC))$ of 
index $4$ that is not an extremal element in $CP_{\sigma}(\IM_3(\IC))$.

\vsp 
\NI (c) For $d=3$, it is obvious that if $\tau$ an extremal element in $CP_{\sigma.\phi_0}(\IM_3(\IC))$ of index $4$ then $\tau$ can not be an extremal in $CP_{\sigma}(\IM_3(\IC))$ 
as $16 > 9$. H. Ohno [Par,Oh] gave an upper limit for $d$ and in case of $n=4$, $d$ can be atmost $4$. Now we aim to classify all extremal elements $\tau$ in $CP_{\sigma,\phi_0}(\IM_3(\IC))$ upto coclycle conjugacy.   

\vsp 
Let $\tau$ be an extremal element in $CP_{\sigma,\phi_0}(\IM_3(\IC))$ and $\tau = \lambda \tau_1 + (1-\lambda)\tau_0$ be its decomposition into extremal elements in $CP_{\sigma}(\IM_3(\IC)$, where $\tau_0$ and $\tau_1$ are two extremal element in $CP_{\sigma}(\IM_3(\IC)$ of index $2$ such that linear span of $\clm_{\tau_0}$ and $\clm_{\tau_1}$ is $4$. For a proof, we can employ argument used in (b). 

\vsp 
Now we will prove the converse statement. Let $\tau_0$ and $\tau_1$ be two extremal elements in $CP_{\sigma}(\IM_3(\IC))$ of index $2$ but not elements in $CP_{\sigma,\phi_0}$ such that $\clm_{\tau_0}$ and $\clm_{\tau_1}$ span $4$ dimensional vector space $\clm_{\tau}$, where $\tau=\lambda \tau_1 + (1-\lambda)\tau_0$ is an element in $CP_{\sigma,\phi_0}$ for some $\lambda 
\in (0,1)$. We claim that $\tau$ is an extremal element in $CP_{\sigma,\phi_0}(\IM_3(\IC))$. 

\vsp 
We may write for $x \in \IM_3(\IC)$ that 
$$\tau(x)=v_1xv_1^*+v_2xv_2^*+v_3xv_3^*+v_4xv_4^*,$$ 
where $v^*_1,v^*_2$ spans $\clm_{\tau_0}$, $v_3^*,v_4^*$ spans $\clm_{\tau_0}$, 
$v_1v_1^*+v_2v_2^*=\lambda I_3$ and $v_3v_3^*+v_4v_4^* =(1-\lambda)I_3$ for some $\lambda \in (0,1)$. Following (1.11.3) in [Ar2], we may choose such a family $(v_i)$ and an inner-product $s$ on the vector space $\clm_{\tau}$ with $s(v^*_i,v_j^*)=\delta^i_j$ for $1 \le i,j \le 4$ and $\phi_0(v_i^*v_j)=0$ for $i \neq j$. 

\vsp 
Suppose $\tau$ is not an extremal element in $CP_{\sigma}(\IM_3(\IC))$. Then we may write $\tau =\sum_{0 \le k \le m}\mu_k \eta_k$ for some extremal elements $\eta_k \in CP_{\sigma,\phi_0}(\clm)$, where $\mu_k \in (0,1)$ and $\sum_{0 \le k \le m} \mu_k=1$. By Proposition 4.2, index of $\eta_k$ is $4$ for some $0 \le k \le m$ if and only if $\eta=\eta_k$ for all $0 \le k \le m$. So index of $\eta_k \le 3$ for each $0 \le k \le m$. 

\vsp 
If index of $\eta_0$ is $3$ then for $x \in \clm$, $\mu_0 \eta_0(x)=\sum_{1 \le k \le 3} w_kxw_k^*$ and index of $\tau_1$ is $1$ where $s(w^*_i,w_j^*)=\delta^i_j$ for all $1 \le i,j \le m$, $\phi_0(w_iw_j^*)=0$ for all $i \ne j$ and $u_4=c^{-1}_4w_4$ is a unitary matrix for $|c_4|^2=\mu_1=1-\mu_0$. Recall $\mu_1$ is taken to be the maximal values for $c > 0$ for which the map $x \raro \tau(x)-cu_4xu_4^*$ is completely positive. 

\vsp 
So $w^*_4 \in \clm_{\tau}$ but $w^*_4 \notin \clm_{\tau_0}$ and $w_4^* \notin \clm_{\tau_1}$. Without loss of gernerality, we may assume that the linear span of $\{v^*_1,v^*_4\}$ is equal to the linear span of $\{ w^*_1,w^*_4 \}$ if needed by changing the orthonomal bases of $\clm_{\tau_0}$ and $\clm_{\tau_1}$. If so then the linear span of $\{w^*_2,w^*_3\}$ is equal to the 
linear span of $\{v^*_2,v^*_3\}$.

\vsp 
So the set $\{w^*_k:1 \le k \le 4\}$ is also an orthonormal basis for $\clm_{\tau}$ in the 
inner product $s$. In particular, we have 
\be 
v^*_j = \sum_{1 \le k \le 4}a^j_kw_k^*
\ee 
for some $A =[a^j_k] \in U_4(\IC)$, where 
$$A= \left ( \begin{array}{llll}&a^1_1,& \;\;\;0, \;\;\;0,\;\;\;a^1_4 \;\; \\ &0,&\;\;a^2_2,\;\;a^2_3, \;\;\; 0 \;\;\;\\ &0,&\;\;a^3_2,\;\;\;a^3_3,\;\;\;0\;\;\;\\ 
&a^4_1,&\;\;\;0,\;\;\;0,\;\;\;\;a^4_4\;\;\; \\
\end{array} \right )
$$

So we also have
\be 
v_j =\sum_{1 \le k \le 4}\bar{a}^j_kw_k
\ee
The unital CP map $\tilde{\tau}(x)=\sum_k v^*_kxv_k$ is also equal to $\mu_0 \tilde{\eta}+\mu_1\tilde{\eta}_1$, where $\tilde{\eta}_i$ are dual of $\eta_i$ for $0 \le i \le 1$. The element $\tilde{\tau}$ can not be an extremal element and it admits a decomposition 
i.e. $\tilde{\tau}(x)=\sum_{1 \le k \le 4}l_kxl_k$, where $l_1l_1^*+l_2l_2^*=\nu I$ and $l_3l_3^*+l_4l_4^*=1-\nu$. Now we repeat the argument used for $\tau$ and assume without loss of generality that the linear spans of $\{l^*_1,l^*_4 \}$ and $\{l^*_2,l^*_3\}$ are equal to the linear spans of $\{w_1,w_4\}$ and $\{w_2,w_3\}$ respectively. So the linear span of 
orthogonal elements $\{l_1^*,l_2^*\}$ is equal to linear span of orthogonal elements 
$\{v_1,v_2 \}$ by (24). In particular, $\sum_{1 \le k \le 2}v_k^*v_k = \sum_{1 \le k \le 2}l_kl_k^*=\nu I$. Now we use the tracial property of $\phi_0$ to compute 
$\nu_0 = \sum_{1 \le k \le 2}\phi_0(v_k^*v_k) = \sum_{1 \le k \le 2}\phi_0(v_kv_k^*)=\lambda$. This brings a contradition to our hypothesis that $\tau_0$ and $\tau_1$ are not elements in $CP_{\sigma,\phi_0}$.

\vsp 
We sum up our results on extreme points of $CP_{\sigma,\phi_0}(\IM_3(\IC))$ 
in the following theorem.

\vsp 
\begin{thm} 
Let $\tau$ be an element in $CP_{\sigma,\phi_0}(\IM_3(\IC))$ of index $d$. Then $1 \le d \le 4$ and the following statements hold:

\vsp 
\NI (a) For $1 \le d \le 3$, $\tau$ is an extremal element in $CP_{\sigma,\phi_0}(\IM_3(\IC))$ if and only if $\tau$ is also an extremal element in $CP_{\sigma}(\IM_3(\IC))$;

\vsp 
\NI (b) For $d=4$, $\tau$ is an extremal element in $CP_{\sigma,\phi_0}(\IM_3(\IC))$ if and only if $\tau = \lambda \tau_1 + (1-\lambda)\tau_0$ for some $\lambda \in (0,1)$, where $\tau_0$ and $\tau_1$ are extremal elements in $CP_{\sigma}(\IM_3(\IC))$ but are not elements in $CP_{\sigma,\phi_0}(\IM_3(\IC))$. 
\end{thm} 

\end{exam}

\vsp 
\begin{rem} 
Theorem 4.6 (b) gave not only a classification of extremal elements in $CP_{\sigma,\phi_0}(\IM_3(\IC))$ of index $4$ but also a constructive method for an extremal element in $CP_{\sigma,\phi_0}(\IM_3(\IC))$ of index $4$, thus not an extremal in $CP_{\sigma}(\IM_3(\IC))$. Same method can be employed to classify extremal elements for $CP_{\sigma,\phi_0}(\IM_3(\IC))$ in an inductive scheme.  
\end{rem}

\bigskip
{\centerline {\bf REFERENCES}}

\begin{itemize} 

\bigskip
\item{[Ar1]} Arveson, W.: Sub-algebras of $C^*$-algebras, Acta Math. 123, 141-224, 1969. 

\item {[Ar2]} Arveson, W.: On the index and dilations of completely positive semigroups. Internat. J. Math. 10 (1999), no. 7, 791--823. 

\item {[Bi]} Birkhoff, G. : Three observations of linear algebra, Rev.Uni.Nac. Tucuman (A), 5(146):147-151.  

\item {[BR]} Bratteli, Ola., Robinson, D.W. : Operator algebras and quantum statistical mechanics, I,II, Springer 1981.

\item {[Ch]} Choi, M.D.: Completely positive linear maps on complex matrices, Linear Algebra and App (10) 1975 285-290. 

\item{[ChC]} Choi, Man Duen, Christensen, E.:
Completely order isomorphic and close $C^*$-algebras need not be $*$-isomorphic.
Bull. London Math. Soc. 15 (1983), no. 6, 604-610.

\item {[ChE]} Christensen,E., Evans,D.: Cohomology of operator algebras and quantum dynamical semigroups, J. London Math. Soc. 20 (1979), 358--368. 

\item{[Da]} Davies, E.B.: Quantum Theory of open systems, Academic press, 1976.

\item{[ER]}  Effros, Edward G.; Ruan, Zhong-Jin: Operator spaces. London Mathematical Society Monographs. New Series, 23. The Clarendon Press, Oxford University Press, New York, 2000. 

\item {[EvL]} Evans, D., Lewis, J. T.: Dilations of irreversible evolutions in algebraic quantum theory, Comm. Dubl. Inst. Adv. Studies, Ser A24 (1977). 

\item {[Fa]} Farenick, Douglas; Arveson's criterion for unitary similarity. Linear Algebra Appl. 435, no. 4, 769-777 (2011).

\item {[Ha1]} Hamana, M.: Injective envelopes of $C^*$ -algebras. J. Math. Soc. Japan, 31, 181-197, (1979). 

\item {[Ha2]} Hamana, M.: Injective envelopes of operator systems. Publ. Res. Inst. Math. Sci. (Kyoto), 15, 773-785, (1979).  

\item {[Ka1]} Kadison, Richard V.: Isometries of operator algebras, Ann. Math. 54(2)(1951) 325-338.  

\item {[Ka2]} Kadison, Richard V.: A generalized Schwarz inequality and algebraic invariants for operator algebras,  Ann. of Math. (2)  56, 494-503 (1952). 

\item {[Ke]} Kendall, D.G. : On infinite doubly stochastic matrices and Birkhoff's problem, III. London Math. Soc. J., 35 (1960):81-84. 

\item {[Ko]} K\"{o}nig, D. : The theory of finite and infinite graphs, T\"{a}ubner 1936, Birkh\"{a}ser, Boston, 1990, p-327. 

\item {[KuM]} Kümmerer, B., Maassen, H. : The essentially commutative dilations of dynamical semigroups on $M_n$. Comm. Math. Phys. 109 (1987), no. 1, 1--22. 

\item {[LS]} Landau, L.J., Streater, R.F.: On Birkhoff theorem for doubly stochastic completely positive maps of matrix algebras, 
Linear Algebra and its Applications, Vol 193, 1993, 107-127 

\item{[LeL]} Leung, Denny H. ; Li, Lei: Order isomorphisms on function spaces. Studia Maths. 219 (2013), no. 2, 123-138. 

\item {[MW]} Mendl, Christian B., Wolf, Michael M.: Unital quantum channel's convex structure and revivals of Birkhoff's theorem. Comm. Math. Phys. 289 (2009), no. 3, 1057--1086. 

\item {[Mo1]} Mohari, A., Pure inductive limit state and Kolmogorov's property. II, J. Operator Theory 72 (2014), no. 2, 387-404. 

\item {[Mo2]} Mohari, A., A mean ergodic theorem of an amenable group action,
Infin. Dimens. Anal. Quantum Probab. Relat. Top. 17 (2014), no. 1, 1450003, 13 pp.

\item {[Mo3]} Mohari, A., Extremal unital completely positive maps and their symmetries,
Complex Anal. Oper. Theory 12 (2018), no. 7, 1739-1765.

\item{[Oh]} Ohno, H.: Maximal rank of extremal marginal tracial states.  J. Math. Phys.  51  (2010),  no. 9, 092101, 9 pp. 

\item {[Pa]} Paulsen, V.: Completely bounded maps and operator algebras, Cambridge Studies in Advance Mathematics 78, 
Cambridge University Press. 2002 

\item{[Par]} Parthasarathy, K.R.: extremal quantum states in coupled states in coupled systems, Ann.Inst. H. Poincar\'{e}, 41, 257-268 (2005).  
 
\item {[Phe]} Phelps, R.R.: Lectures on Choquet's theorem, vol-1757 of Lecture Notes in Mathematics, Springer-Verlag, Berlin, Second edition (2001). 

\item {[Pi]} Pisier, G.: Introduction to operator space theory. London Mathematical Society Lecture Note Series, 294. Cambridge University Press, Cambridge, 2003. viii+478 pp.

\item{[PSa]} Price, G.L., Sakai, S.: Extremal marginal tracial states in couple systems, Operators and matrices, 1, 153-163 (2007). 

\item{[Ru]} Rudolph, O.: On extremal quantum states of composite systems with fixed marginals, J. Math. Phys. 45, 4035-4041 (2004).

\item {[St]} Stinespring, W. F.: Positive functions on $C^*$ algebras, Proc. Amer. Math. Soc. 6 (1955) 211-216. 
 
\item {[Sto1]} Stone, M. H.: The theory of representations for Boolean algebras. Trans. Amer. Math. Soc. 40 (1936), no. 1, 37-111.

\item {[Sto2]} Stone, M. H.: Applications of the theory of Boolean rings to general topology. Trans. Amer. Math. Soc. 41 (1937), no. 3, 375-481.

\item {[St\o 1]} St\o rmer, E.: Positive linear maps of operator algebras. Acta Math. 110, 233-278 (1963).

\item{[St\o 2]} St\o rmer, E.: Positive linear maps of operator algebras. Springer Monographs in Mathematics. Springer, Heidelberg, 2013.

\end{itemize}

\end{document}